\documentclass[preprint,sort&compress]{elsarticle}
\usepackage{graphicx}
\usepackage{amsmath,amsfonts,amssymb,mathtools,amsthm}
\usepackage{url}
\usepackage{color}
\usepackage{subcaption}
\usepackage{algorithm}
\usepackage[noEnd=true]{algpseudocodex}
\usepackage{float}
\usepackage[shortlabels]{enumitem}
\usepackage{multirow}
\usepackage{diagbox}

\usepackage{xr-hyper}
\usepackage[hidelinks]{hyperref}
\usepackage{cleveref}
\usepackage{todonotes}

\def\tto{\;{\lower 1pt \hbox{$\rightarrow$}}\kern -10pt
	\hbox{\raise 2pt \hbox{$\rightarrow$}}\;}

\def\Bar{\overline}
\def\ra{\rangle}
\def\la{\langle}

\def\R{\mathbb{R}}
\def\N{I\!\!N}

\def\ox{\bar{x}}

\def\bX{\mathbf{X}}
\def\bO{\mathbf{\Omega}}
\def\bA{\mathbf{A}}
\def\bE{\mathbf{E}}

\def\dom{\mbox{\rm dom}\,}

\def\O{\Omega}
\def\ph{\varphi}

\newcommand\norm[1]{\left\|#1\right\|}

\def\emp{\emptyset}

\def\oR{\Bar{\R}}

\def\ri{\text{ri}}

\newtheorem{proposition}{Proposition}
\newtheorem{theorem}{Theorem}
\newtheorem{corollary}{Corollary}

\theoremstyle{definition}
\newtheorem{Example}{Example}

\numberwithin{equation}{section} 

\title{The Boosted DC Algorithm for Clustering with Constraints}

\author[1]{Tuyen Tran \corref{cor1}}
\ead{ttran18@luc.edu}

\author[2]{Kate Figenschou}
\ead{katefig@googlemail.com}

\author[2]{Phan Tu Vuong}
\ead{t.v.phan@soton.ac.uk}

\affiliation[1]{organization={Department of Mathematics and Statistics, Loyola University Chicago}, 
  postcode={60660},
  city={Chicago},
  country={U.S.A}
}

\affiliation[2]{organization={Mathematical Sciences, University of Southampton},
  postcode={SO17 1BJ},
  city={Southampton},
  country={UK}
}

\cortext[cor1]{Corresponding Author}

\makeatletter
\def\ps@pprintTitle{%
  \let\@oddhead\@empty
  \let\@evenhead\@empty
  \def\@oddfoot{\reset@font\hfil\thepage\hfil}
  \let\@evenfoot\@oddfoot
}
\makeatother

\begin{document}
	\begin{abstract} This paper aims to investigate the effectiveness of the recently proposed Boosted Difference of Convex functions Algorithm (BDCA) when applied to clustering with constraints and set clustering with constraints problems. This is the first paper to apply BDCA to a problem with nonlinear constraints. We present the mathematical basis for the BDCA and Difference of Convex functions Algorithm (DCA), along with a penalty method based on distance functions. We then develop algorithms for solving these problems and computationally implement them, with publicly available implementations. We compare old examples and provide new experiments to test the algorithms. We find that the BDCA method converges in fewer iterations than the corresponding DCA-based method. In addition, BDCA yields faster CPU running-times in all tested problems.
	\end{abstract}
	

	\maketitle
	\section{Introduction}
	In the field of mathematical optimization, the properties of convex functions and convex sets have allowed the development of numerical algorithms to solve problems efficiently. Generally, the aim of an optimization problem is to minimize an objective function, with respect to some constraints in order to find the best possible (and hence smallest) objective values. A local minima will be at least as good as any nearby elements, while a global minima will be at least as good as every feasible element. Many objective functions have several local minima, which makes identifying global minima difficult, but the properties of convex functions are particularly useful in this context. For a convex function, if there is a local minimum that is interior, it is also the global minimum.
	
	However, while convex functions can be useful modeling tools, most real-world problems are non-convex. Such problems are generally more complicated and difficult. Indeed, most non-convex optimization problems are at least NP-hard. Non-convexity presents many challenges but in particular the presence of both local and global minima, and the lack of identifiable characteristics for global minima greatly increases the computational complexity. Various methods have been developed to tackle these types of problems, which can broadly be split into either global or local approaches. Global approaches, such as branch and bound, are very expensive (especially for large-scale problems) but are able to guarantee the globality of the solution. Local approaches meanwhile are faster and cheaper, but their solutions cannot generally be proven to be global. Even local approaches struggle to be effective at large-scale too, so the challenge of developing algorithms which balance quality and scalability is complex.
	
	DC programming is a non-convex optimization problem where the objective function is a Difference of Convex (or DC) function. The approach uses the convexity of the DC components and duality to make solving the nonconvex problem easier. It has been shown to be robust and efficient in many applications, including with large-scale problems, and is relatively simple to use and implement \cite{An2005,TA1}. In the last 30 years, the DC Algorithm (DCA) is the golden method for solving DC programming.  The computed solutions cannot be guaranteed to be global as the DCA converges to a local solution, however in experiments the DCA often converges to a global solution \cite{TA1}. The method was initially introduced by Pham Dinh Tao in 1985, following on from his work on subgradient algorithms for convex maximization programming, and further developed by Le Thi Hoai An and Pham Dinh Tao in the 1990's \cite{LeThi2018}. In the following years, the DCA was applied to various different topics, especially in machine learning and data mining problems, becoming increasingly popular.
	
	Following on from the success and popularity of DC programming, new algorithms based on the DCA have been proposed. The Boosted Difference of Convex functions Algorithm (BDCA) is one of these new methods, introduced to accelerate the convergence of the classical DCA \cite{AragonArtacho2018}. More importantly, BDCA escapes from bad local solutions thanks to the line search step with arbitrary large trial stepsize. Numerical experiments have shown that the BDCA is able to outperform the DCA in problems such as Minimum Sum-of-Squares Clustering (MSSC) and $l_\infty -$trust-region subproblems \cite{Artacho2020}. Whether the BDCA can be successfully applied in other settings is a topic of ongoing research, and the aim of this project was to investigate the use of the BDCA when applied to clustering problems with constraints. \(\)
	
	Clustering is a common statistical data analysis method which aims to classify similar objects together into groups (or clusters). There are many different approaches to defining the similarity of objects and how they are assigned to clusters, and generally no single algorithm will be correct for a given task. Centroid-based algorithms represent clusters with a central vector, assigning clusters based on proximity to cluster centers using a proximity metric. The k-means algorithm, a centroid-based approach with a fixed number of $k$ clusters, uses the squared Euclidean distance and is one of the most widely known methods. However, though it is simple and easy to implement, it suffers from certain weaknesses including a high dependence on the initial choice of cluster centers as well as the proximity measure, and the algorithm has no guarantee of convergence to a global optimal solution. Much research has focused on alternative algorithms to k-means which alleviate its drawbacks, including DC programming based approaches.
	
	\subsection*{The development of DC programing and DCA}
	Pham Dinh Tao and Le Thi Hoai An, who have extensively developed DC programming and the DCA, have published many papers on the topic. Their 1997 paper "DC programming and DCA: Theory, Algorithm, Applications" presented the most complete study of DC programming and the DCA at that point \cite{TA1}. It described key components of topic including DC duality local optimality conditions, convergence properties of the DCA and its theoretical basis. Alongside the new and significant mathematical results, the paper presented numerical experiments applied to real-life problems which proved the effectiveness of the DCA compared to known algorithms. The extensive content of the paper made it a highly important source for further work on the topic.
	
	Following on from the development of DC programming in the late $20^{\text{th}}$ century, the early 2000's saw increasing applications of the DCA, especially in machine learning and related fields. The first paper to investigate the use of the DCA for clustering focused on a K-median clustering problem and the K-means algorithm \cite{An2007}. When tested with real-world databases, both of the DCA methods presented in the paper achieved better objective values than the classic clustering algorithms and were faster, too. The success and efficiency of the methods showed that clustering was an area of interest for further research into the DCA and led to further papers looking into different clustering types, including minimum sum-of-squares clustering (MSSC), fuzzy clustering and hierarchical clustering.
	
	The following years saw DC programming and the DCA become classical, with applications across a wide variety of fields. Further investigation into the use of the DCA to solve clustering problems continued. The clustering and multifacility location problems with constraints described in this paper were originally developed in \cite{Nam2018}. This paper presented the mathematical basis for the DC decompositions of each model, including the use of a penalty method based on the squared distance function. Numerical examples were presented for the models, though extensive testing on the effectiveness of the methods was not included.
	
	As the DCA become increasingly popular, research was performed on methods to accelerate its performance and recently the BDCA was proposed. The addition of a line search step in BDCA was proven to accelerate the convergence of DCA. Moreover, it has been noticed that BDCA escapes from bad local solutions, which cannot be the case in DCA.
	More details about increasing performance, even in high-dimensions can be found in \cite{AragonArtacho2018}. Further developments on the use of the BDCA for linearly constrained problems \cite{AragonArtacho2022} and nonsmooth functions \cite{Artacho2020} have continued to show the effectiveness of the BDCA in various applications. The BDCA has been applied to Minimum Sum-of-Squares Clustering, where it was on average 16 times faster than the DCA \cite{Artacho2020}. Further extensions of BDCA can be found in  \cite{Ferreira2021,Ferreira2022}. The aim of this paper is to explore whether the BDCA is also effective for constrained clustering problems.
	
	In this paper, we investigate the performance of the BDCA against the DCA for solving clustering problems with constraints and set clustering problems. The DCA for solving these problems were originally developed in \cite{Nam2018}, and this project follows from the work there. This work is split between studying the mathematical basis for the algorithms presented, and testing implementations of these algorithms in MATLAB.
	
	The paper is structured as follows.  Firstly, \cref{Chapter:Prelim} presents some basic mathematical tools of convex analysis, and followed by explanation of the DCA and BDCA. Next, the penalty method is reviewed in \cref{Chapter:Penalty}. \Cref{Chapter:Clusterings} lays out the first clustering with constraints problem, while in \cref{Chapter:SetClusters} we study a model of clustering with constraints involving sets. The numerical tests performed are explained in \cref{examples}. Finally, \cref{Chapter: Conclusions} summarizes the findings of the paper as well as future research directions.
	
	\section{Preliminaries}\label{Chapter:Prelim}
	In this section, we present basic tools of analysis and optimization. The readers are referred to \cite{An2007,HUL,bn,r} for more details and proofs of the presented results.
	
	Let us define $\oR \coloneqq \R \cup \{+\infty\}=(-\infty, \infty]$ and let $f\colon \R^d\to \oR$ be a convex function. An element $v\in\R^d$ is called a \emph{subgradient} of $f$ at $\bar x\in \dom(f)=\{x\in \R^d\; |\; f(x)<\infty\}$ if it satisfies
	\begin{equation}
		\langle v,x-\bar x\rangle\leq f(x)-f(\bar x)\; \mbox{\rm for all } x\in \R^d.
	\end{equation}
	The set of all such elements $v$ is called the \emph{subdifferential} of $f$ at $\bar x$ and is denoted by $\partial f(\bar x)$. If $\ox\not\in \dom(f)$, we set $\partial f(\ox)=\emptyset$. Subdifferentials possess many calculus rules that are important in practice. In particular, for a finite number of convex functions  $f_i\colon\R^d\to \oR$, $i=1, \ldots, m$, we have the following sum rule:
	\begin{equation}
		\partial(f_1+\cdots+f_m)(\ox)=\partial f_1(\ox)+\cdots+\partial f_m(\ox)\; \mbox{\rm for all }\ox\in \R^d
	\end{equation}
	provided that $\bigcap _{i=1}^m \ri (\dom(f_i)) \neq \emptyset$, where $\ri(\O)$ denotes for the \emph{relative interior} of $\O$; see, e.g, \cite[Definition~1.68]{bn}.
	
	If $f=\max_{i=1, \ldots, m} f_i$,  and $f_i$ is continuous at $\ox$ for every $i=1,\ldots,m$, then for any $\ox\in \R^d$ we have the following maximal rule:
	\begin{equation}\label{convexmax}
		\partial f(\ox)=\mbox{\rm conv}\left( \bigcup_{i\in I(\ox)} \partial f_i(\ox)\right),
	\end{equation}
	where $I(\ox)=\{i\;|\; f_i(\ox)=f(\ox)\}$, and \(conv \) is the \emph{convex hull}.
	
	Given a nonempty closed convex subset $\O$ of $\R^d$ with $\ox\in \O$, the \emph{normal cone} to $\O$ at $\ox$ is defined by
	\begin{equation}
		N(\ox; \O)=\big\{v\in \R^d\;\big|\; \la v, x-\ox\ra \leq 0\; \mbox{\rm for all }x\in \O\big\}.
	\end{equation}
	
	If $\ox\not\in \O$, we set $N(\ox,\O)=\emptyset$. It is well-known that an element $\ox\in \R^d$ is an absolute minimizer of a convex function $f\colon \R^d\to \R$ on $\O$ if and only if $\ox$ is a local minimizer of $f$ on $\O$. Moreover, this happens if and only if the following optimality condition holds:
	\[
	0\in \partial f(\ox) +N(\ox; \O).
	\]
	Let $\Theta\subset\R^d$ be a nonempty set (not necessarily convex). The \emph{distance function} to $\Theta$ is defined by
	\begin{equation*}
		d(x; \Theta)=\inf\left\{\|x-w\|\; \big|\; w\in \Theta\; x\in \R^d\right\}.
	\end{equation*}
	The \emph{Euclidean projection} from $x\in \R^d$ to $\Theta$ is the set
	\begin{equation*}
		P(x; \Theta)=\big\{w\in \Theta\; \big|\; d(x; \Theta)=\|x-w\|\big\}.
	\end{equation*}
	There are two important properties of the Euclidean projection. First, if $\Theta$ is a nonempty closed set, then $P(x; \Theta)$ is nonempty and is a singleton if $\Theta$ is also convex. Second, if $\Theta$ is a convex set and $w\in P(x; \Theta)$, then $x-w\in N(w; \Theta)$.
	
	Another tool we will use is the notion of \emph{Fenchel conjugates}. Let $g\colon \R^d\rightarrow\R$ be a function (not necessarily convex). The Fenchel conjugate of $g$ is defined by
	\begin{align*}
		g^*(y)=\sup\left\{\langle y,x\rangle-g(x)\; \big|\; x\in\R^d\right\},\; y\in \R^d.
	\end{align*}
	Note that $g^*\colon\R^d\to \oR$ is an extended-real-valued convex function. Suppose further that $g$ is itself convex, then the {\em Felchel-Moreau theorem} states that $(g^*)^*=g$. Based on this theorem, we have the following relation between the subgradients of $g$ and its Fenchel conjugate:
	\begin{equation}\label{fcg}
		x\in\partial g^*(y)\iff y\in\partial g(x).
	\end{equation}
	
	The notions of subgradient and Fenchel conjugate provide the mathematical foundation for the DCA introduced in the next section. The following proposition gives us a two-sided relationship between the Fenchel conjugates and subgradients of convex functions.
	
	\begin{proposition}\label{characterization} Let $\ph\colon\R^d\to\oR$ be a proper, lower semicontinuous, and convex function. Then $v\in\partial\ph^*(y)$ if and only if
		\begin{equation*}\label{e1}
			v\in\mbox{\rm argmin}\,\big\{\ph(x)-\la y,x\ra\;\big|\;x\in\R^d\big\}.
		\end{equation*}
		Furthermore, $w\in\partial\ph(x)$ if and only if
		\begin{equation*}\label{e2}
			w\in\mbox{\rm argmin}\,\big\{\ph^*(y)-\la x,y\ra\;\big|\;y\in\R^d\big\}.
		\end{equation*}
	\end{proposition}
	
	The proof of this proposition can be found in \cite[Proposition~2.1]{bmnt19}.
	
	Throughout this paper, we denote $\bA\in\R^{m\times d}$  as the \emph{data matrix}. The $i^{th}$ row is denoted $a^i\in \R^d$ for $i=1,\dots,m$. Similarly, $\bX\in\R^{k\times d}$ is defined as the \emph{variable matrix}  and the $\ell^{th}$ row is denoted $x^\ell\in \R^d$ for $\ell=1,\dots,k$.
	The linear space $\R^{k\times d}$ is equipped with the inner product $\langle \bX, \mathbf Y\rangle=\mbox{\rm trace}(\bX^T \mathbf Y)$.
	
	Recall that the \emph{Frobenius norm} on $\R^{k\times d}$ is defined by
	\begin{align*}
		\big\| \bX\big\|_F= \left\langle \bX,\bX\right\rangle^{1/2}=  \left(\sum_{\ell =1}^{k}\langle x^\ell, x^\ell\rangle\right)^{1/2}= \left(\sum_{\ell =1}^{k}\|x^\ell\|^2\right)^{1/2}.
	\end{align*}
	Notice that the squared Frobenius norm is differentiable with the following representation
	\begin{equation*}
		\nabla\norm{\bX}^2_F=2\bX\; \mbox{\rm for }\bX\in\R^{k\times d}.
	\end{equation*}
	
	For constraint sets \(\O^\ell\) used in this paper, we are using the same notation and assumptions as those sets as in \cite{Nam2018}. In which, $\O^\ell \subset \R^d$ for $l=1,\ldots, k$ are nonempty closed convex sets and the Cartesian set product is defined as $\bO=\O^1\times \O^2\times \ldots\times \O^k.$ For $\bX \in \R^{k\times d}$, the projection from $\bX$ to $\bO$ is the matrix $\mathbf Y$ whose $\ell^{th}$ row is $y^\ell=P(x^\ell; \O^\ell).$
	The relationship between distance function and Frobenious norm can be viewed as
	\begin{equation*} \label{disXO}
		[d(\bX; \bO)]^2=\|\bX - \mathbf Y\|_F^2= \sum_{\ell =1}^{k}\|x^\ell - y^\ell\|^2=\sum_{\ell =1}^{k} d(x^\ell; \O^\ell)^2.
	\end{equation*}
	\section{DCA and BDCA}
	\subsection{The DCA}\label{Prelim:DCA}
	
	The notions of subgradients and Fenchel conjugates in the previous section provides  mathematical foundation for the DCA introduced below.
	Consider the difference of two convex functions $g-h$ on a finite-dimensional space and assume that $g\colon\R^d\to\oR$ is extended-real-valued while $h\colon\R^d\to\R$ is real-valued on $\R^d$. Then a general problem of {\em DC optimization} is defined by:
	\begin{equation}\label{DCP}
		\mbox{\rm minimize }f(x):=g(x)-h(x),\quad x\in\R^d. \tag{\ensuremath{\mathcal{P}}}
	\end{equation}
	
	The DCA introduced by Tao and An is a simple but effective algorithm for minimizing the function $f$; see \cite{TA1,TA2}.
	
	
	\begin{algorithm}[H]
		\caption{N Iteration DCA} \label{alg1}
		\begin{algorithmic}[0]
			\Procedure{DCA}{$x_1, N \in \mathbb{N}$}
			\For{\(p = 1, \dots , N\)}
			\State Find $y_p \in \partial h(x_p)$
			\State Find $x_{p+1} \in \partial g^*(y_p)$
			\EndFor
			\EndProcedure
			\State \Output{$x_{N+1}$}
		\end{algorithmic}
	\end{algorithm}
	
	To proceed further, recall that a function $\ph\colon\R^d\to\oR$ is $\gamma$-\emph{convex} with a given modulus $\gamma\ge 0$ if the function $\psi(x):=\ph(x)-\frac{\gamma}{2}\|x\|^2$ as $x\in\R^d$ is convex on $\R^d$. If there exists $\gamma>0$ such that $\ph$ is $\gamma-$convex, then $\ph$ is called \emph{strongly convex} on $\R^d$.
	
	We also recall that a vector $\ox\in\R^d$ is a \emph{stationary point}/\emph{critical point} of the DC function $f$ from \cref{DCP} if
	\begin{equation*}
		\partial g(\ox)\cap\partial h(\ox)\ne\emp.
	\end{equation*}
	The next result, which can be derived from \cite{TA1,TA2}, summarizes some convergence results of the DCA. Deeper studies of the convergence of this algorithm and its generalizations involving the  Kurdyka-Lojasiewicz (KL) inequality are given in \cite{AnNam,TAN}.
	
	\begin{theorem}\label{dsa-conv} Let $f$ be a DC function taken from \eqref{DCP}, and let $\{x_k\}$ be an iterative sequence generated by Algorithm~{\rm 1}. The following assertions hold:
		\begin{enumerate}
			\item The sequence $\{f(x_k)\}$ is always monotone decreasing.
			\item Suppose that $f$ is bounded from below, that $g$ is lower semicontinuous and $\gamma_1$-convex, and that $h$ is $\gamma_2$-convex with $\gamma_1+\gamma_2>0$. If $\{x_k\}$ is bounded, then the
			limit of any convergent subsequence of $\{x_k\}$ is a stationary point of $f$.
		\end{enumerate}
	\end{theorem}
	
	In many practical applications of Algorithm~1, for a given DC decomposition of $f$ it is possible to find subgradient vectors from $\partial h(x_k)$ based on available formulas and calculus rules of convex analysis. However, it may not be possible to explicitly calculate an element of $\partial g^*(y_k)$. Such a situation requires either constructing a more suitable DC decomposition of $f$, or finding $x_{k+1}\in\partial g^*(y_k)$ approximately by using the minimization criterion of \cref{characterization}. This leads us to the following modified version of the DCA.
	
	\begin{algorithm}[H]
		\caption{N Iteration for \bf{DCA-2}} \label{algDCA2}
		\begin{algorithmic}
			\Procedure{DCA-2}{$x_1 \in \mbox{\rm dom}\,g, N \in \mathbb{N}$}
			\For{p = 1, \dots , N}
			\State Find $y_p \in \partial h(x_p)$
			\State Find $x_{p+1}$ approximately by solving the problem 
			\[
			\mbox{\rm minimize}\;\ph_p(x):=g(x)-\la y_p,x\ra,\;x\in\R^d.
			\]
			\EndFor
			\EndProcedure
			\State \Output{$x_{N+1}$}
		\end{algorithmic}
	\end{algorithm}
	
	\subsection{The BDCA}\label{Prelim:BDCA}
	
	The Boosted DC Algorithm (BDCA) has been recently proposed to accelerate the performance of the DCA. The BDCA has an extra line search step at the point found by the DCA at each iteration. This allows the BDCA to take larger steps leading to a larger reduction of the objective value each iteration. The BDCA has also been found to escape bad local optima more easily than the DCA, leading to better objective values as well as increased speed. It has been proved in the case where $g$ and $h$ are differentiable \cite{AragonArtacho2018}, and when $g$ is differentiable but $h$ is not \cite{Artacho2020}. In this section, the problem \cref{DCP} and associated assumptions when applying the BDCA are presented. The problems in this paper fall under the latter case. There are two assumptions made when applying the BDCA:
	
	\textbf{Assumption 1:} \textit{Both $g$ and $h$ are strongly convex with modulus $\rho > 0$.}
	
	\textbf{Assumption 2:} \textit{The function $h$ is subdifferentiable at every point in $\operatorname{dom} h$. So $\partial h(x) \neq \emptyset$ for all $x \in \operatorname{dom} h$. The function $g$ is continuously differentiable on an open set containing $\dom(h)$ and $\underset{x \, \in \, \mathbb{R}^d}{\operatorname{inf}} \ f (x) > - \, \infty \,$.}\\
	Under these two assumptions, the following optimality condition holds.
	\begin{theorem}[First-Order Necessary Optimality Condition]
		If $x^* \in \dom(f)$ is an optimal solution of \cref{DCP}, then
		\begin{equation} \label{firstordercondition}
			\partial h(x^*) = \lbrace \nabla g(x^*) \rbrace .
		\end{equation}
	\end{theorem}
	
	The proof of this theorem can be found as Theorem 3 in \cite{Toland1979}. Any point satisfying \cref{firstordercondition} is a stationary point of \cref{DCP}. We say that $\Bar{x}$ is a critical point of \cref{DCP} if
	\begin{equation*}
		\nabla g(\Bar{x}) \in \partial h(\Bar{x}).
	\end{equation*}
	
	Every stationary point $x^*$ is a critical point, but in general the converse is not true.
	
	A general form of the BDCA is presented below.
	
	\begin{algorithm}[H]
		\caption{BDCA for solving \cref{P}} \label{AlgorithmGeneralBDCA}
		\begin{algorithmic}
			\Procedure{BDCA}{$x_1$, $\alpha > 0$, $\beta \in (0, 1)$, $N \in \mathbb{N}$}
			\For{$k = 1, \dots , N$}
			\State \underline{Step 1:} Select $u_k \in \partial h(x_k)$ and solve the strongly convex optimization problem 
			\[
			\underset{x \, \in \, \mathbb{R}^d}{\operatorname{min}} \quad \ph_k (x) := g(x) - \langle u_k , x \rangle
			\]
			\State to obtain its unique solution $y_k$.
			\State \underline{Step 2:} Set $d_k := y_k - x_k$
			\If{$d_k = 0$}
			\State \Return $x_k$
			\Else
			\State Go to \underline{Step 3}
			\EndIf
			\State \underline{Step 3:} Choose any $\Bar{\lambda}_k \geq 0$, set $\lambda_k =\Bar{\lambda}_k$
			\While{$f(y_k + \lambda_k d_k) > f(y_k) - \alpha \lambda_k^2 \Vert d_k \Vert^2$}
			\State $\lambda_k = \beta \lambda_k$
			\EndWhile
			\State \underline{Step 4:} Set  $x_{k+1} = y_k + \lambda_k d_k$, and $k = k + 1$
			\EndFor
			\EndProcedure
			\State \Output{$x_{N+1}$}
		\end{algorithmic}
	\end{algorithm}
	
	We can see that if $\lambda_k = 0$, then the iterations of the BDCA and the DCA are the same. Therefore, convergence results for the BDCA can apply in particular to the DCA.
	The following proposition shows that $d_k$ is indeed a descent direction for $f$ at $y_k$.
	
	\begin{proposition}\label{prop: BDCADescentDirection}
		For all $k \in \mathbb{N}$, the following hold:
		\begin{enumerate}[(i)]
			\item \textit{$f(y_k) \leq f(x_k) - \rho \Vert d_k \Vert^2\,$}
			\item \textit{$f'(y_k \, ; \, d_k) \leq - \rho \Vert d_k \Vert^2\,$}
			\item \textit{there is some $\delta_k > 0$ such that}
			\begin{equation*}
				f (y_k + \lambda \, d_k ) \leq f (y_k) - \alpha \, \lambda^2 \, \Vert d_k \Vert^2, \quad \text{for all} \ \lambda \in [0,\delta_k],
			\end{equation*}
			\textit{so that the backtracking step 3 of \cref{AlgorithmGeneralBDCA} terminates finitely.}
		\end{enumerate}
	\end{proposition}
	
	The proof of \cref{prop: BDCADescentDirection}  can be found as Proposition 3.1 in \cite{Artacho2020}. With \cref{prop: BDCADescentDirection} it can be shown that the BDCA results in a larger decrease of the objective function than the DCA at each iteration.  The work presented in \cite{Artacho2020} provides the full mathematical background to the BDCA and the improvement of the performance of DCA in relevant applications.
	
	The first convergence result of the iterative sequence generated by BDCA is presented in the following theorem.
	\begin{theorem}\label{prop:innerloop}
		For any $x_{0}\in \mathbb{R}^{m}$, either BDCA
		returns a critical point of \cref{DCP} or it generates
		an infinite sequence such that the following holds.
		\begin{enumerate}
			\item $f(x_{k})$ is monotonically decreasing and convergent to some
			$f(x^*)$.
			\item Any limit point of $\{x_{k}\}$ is a critical point of \cref{DCP}.
			If in addition, $f$ is coercive then there exits a subsequence
			of $\{x_{k}\}$ which converges to a critical point of\cref{DCP}.
			\item $\sum_{k=0}^{+\infty}\|d_{k}\|^{2}<+\infty.$ Further, if there is some $\overline{\lambda}$ such that $\lambda_k\leq\overline{\lambda}$ for all $k$, then $\sum_{k=0}^{+\infty}\|x_{k+1}-x_{k}\|^{2}<+\infty$.\end{enumerate}
	\end{theorem}
	
	More details of the proof and the convergence under the Kurdyka–Lojasiewicz property can be found in \cite[Section 4]{Artacho2020}.
	\section{A Penalty Method via Distance Functions} \label{Chapter:Penalty}
	
	In this section, we review a penalty method that was introduced in \cite{Nam2018} using distance functions for solving constrained optimization problems and then apply them to DC programming. The \emph{quadratic penalty method} was utilized for this technique; see \cite{CZL12,numopt}. Detailed proofs for theorems and propositions below can be found in \cite{Nam2018}.
	
	We first restate the problem of interest here. Let $f\colon \R^d\to \R$ be a function and let $\O_i$ for $i=1, \ldots, q$ be nonempty closed subsets of $\R^d$ with $\bigcap_{i=1}^q\O_i\neq\emptyset$. Consider the optimization problem:
	\begin{equation}\label{P}
		\begin{array}{ll}
			\mbox{\rm min }&f(x)\\
			\mbox{\rm subject to }&x\in \bigcap_{i=1}^q \O_i.
		\end{array}
	\end{equation}
	The above problem can be rewritten as an unconstrained version given by
	\begin{equation}\label{Q}
		\mbox{\rm min }f_\lambda(x)=f(x)+\frac{\lambda}{2}\sum_{i=1}^q[d(x; \O_i)]^2, \; x\in \R^d.
	\end{equation}
	The following theorem provides a relation between optimal solutions of the  problem \cref{P} and problem \cref{Q} which is obtained by a penalty method based on distance functions. Here $\lambda$ is the penalty term. The proof follows \cite[Theorem~17.1]{numopt}.
	
	\begin{theorem}\label{P1}
		Consider  \cref{P} in which $f\colon \R^d\to \R$ is a l.s.c. function.  Suppose that \cref{P} has an optimal solution. If $\lim_{n\to\infty}\lambda_n=\infty$ and $x_n\in\R^d$ is an absolute minimizer of the function $f_{\lambda_n}$ defined in \cref{Q} for all $n\in \N$, then every subsequential limit of $\{x_n\}$ is a solution of \cref{P}.
	\end{theorem}
	
	Now, let $F\colon \R^{k \times d}\to \R$ be a function and let $\O_i^\ell$ for $\ell=1, \ldots, k$ and $i=1, \ldots, q$ be nonempty closed subsets of $\R^d$. We consider the extended version of \cref{P}
	\begin{equation}\label{P2}
		\begin{array}{ll}
			&\mbox{\rm min}    \qquad \qquad      F(x^1, \ldots, x^k)\\
			& \mbox{\rm subject to }\quad x^\ell\in \bigcap_{i=1}^q\O_i^\ell\;, x^\ell\in \R^d \;\mbox{for} \; \ell=1, \ldots ,k.
		\end{array}
	\end{equation}
	The unconstrained version of \cref{P2} is then given by
	\begin{equation}  \label{Flambda}
		\begin{array}{ll}
			\mbox{\rm min } &F_\lambda(x^1, \ldots, x^k)=F(x^1, \ldots, x^k)+\dfrac{\lambda}{2}\sum \limits_{\ell=1}^k\sum \limits_{i=1}^q[d(x^\ell; \O_i^\ell)]^2\\
			&x^\ell\in \R^d \;\mbox{\rm for } \ell=1, \ldots ,k.
		\end{array}
	\end{equation}
	We denote $\bX=(x^1, \ldots, x^k)\in \R^{k \times d}$ and its $\ell^{\rm th}$ row as $x^\ell$ for $\ell=1, \ldots, k$.
	\begin{corollary} Consider \cref{P2} in which $F\colon \R^{k\times d} \to \R$ is a l.s.c. function. Suppose that \cref{P2} has an optimal solution. If $\lim_{n\to\infty}\lambda_n=\infty$ and $X_n=(x^1_n, \ldots, x^k_n)\in \R^{ k\times d}$ is an absolute minimizer of the function $F_{\lambda_n}$, then every subsequential limit of $\{X_n\}$ is a solution of \cref{P2}. \label{cor32}
	\end{corollary}
	
	Next, we recall a known result on DC decompositions of squared distance functions. The proof can be found in \cite[Proposition 5.1]{Nam2017}.
	
	\begin{proposition}\label{dcr} Let $\O$ be a nonempty  closed set in $\R^d$ (not necessarily convex). Define the function
		\begin{equation*}
			\ph_{\O}(x)=\sup\big\{\la 2 x, w\ra -\|w\|^2 \; \big|\;   w\in \O\big\}=2\sup\big\{\la  x, w\ra -\frac{1}{2}\|w\|^2 \; \big |\;    w\in \O\big\}.
		\end{equation*}
		Then we have the following conclusions:\\
		{\rm\bf (i)} The function $\ph_\O$ is always convex. If we assume in addition that $\O$ is convex, then $\ph_\O$ is differentiable with $\nabla \ph_{\O}(x)=2P(x; \O).$\\
		{\rm \bf (ii)} The function $f(x)=[d(x; \O)]^2$ is a DC function with $f(x)=\|x\|^2-\ph_\O(x)$ for all $x\in \R^d$.
	\end{proposition}
	
	Furthermore, the proof of the contents of the following proposition can be found in \cite[Section 3]{Nam2018}. 
	\begin{proposition} \label{prop34}
		Consider \cref{P} where additionally $f$ is a DC function, and all constraint sets are convex sets that satisfy $\bigcap_{i=1}^q \mbox{ri}(\O_i)\neq \emptyset$. Suppose that $\lim_{n\to\infty}\lambda_n=\infty$ and $x_n$ is a critical  point of the DC function $f_{\lambda_n}$ defined in \cref{Q}. Then every subsequential limit of the sequence $\{x_n\}$ is a critical point of \cref{P}.
	\end{proposition}
	A similar development for relating the problems \cref{P2,Flambda} is shown in \cite[Section 3]{Nam2018}.
	\section{Clustering with Constraints} \label{Chapter:Clusterings}
	
	In this section, we review the problem of \emph{clustering with constraints} in \cite{Nam2018}. The squared Euclidean norm is used for measuring distance. The task is to find $k$ centers $x^1,\ldots,x^k\in \R^d$ of $m$ data points $a^1, \ldots, a^m\in \R^d$ and with the constraint that each $x^\ell\in\bigcap_{i=1}^q{\O_i^\ell}$ for some nonempty closed convex set ${\O_i^\ell}\subset\R^d$ with $\ell=1,\ldots,k$. Without loss of generality,  we can assume  that the numbers of constraints is  equal for each center. The problem of interest is given by
	\begin{equation}\label{Pointclustering}
		\begin{array}{ll}
			\mbox{\rm min } &\psi(x^1, \ldots, x^k)=\sum_{i=1}^m \min_{\ell =1, \ldots, k} \|x^\ell-a^i\|^2\\
			\mbox{\rm subject to } &x^\ell\in\bigcap_{j=1}^q \O_j^\ell\; \mbox{\rm for }\ell=1, \ldots, k.\end{array}
	\end{equation}
	As discussed in the previous section, the unconstrained minimization problem is given by
	\begin{equation}\label{PenPointclustering}
		\begin{array}{ll}
			\mbox{\rm min }&f(x^1, \ldots, x^k)= \;  \dfrac{1}{2}\sum_{i=1}^m \min_{\ell =1, \ldots, k}{ \|x^\ell-a^i\|}^2+ \dfrac{\tau}{2} \sum_{\ell=1}^k\sum_{i=1}^q [d(x^\ell;{\O_i^\ell})]^2,\\
			&x^1, \ldots, x^k\in \R^d,
		\end{array}
	\end{equation}
	where $\tau>0$ is a penalty parameter.
	
	Applying \cref{dcr} for any nonempty closed convex set $\O$ in $\R^d$, and using the \emph{minimum-sum principle} we obtain a DC decomposition of $f$ as below
	\begin{align*}
		f(x^1, \ldots, x^k)&=\Big( \frac{1}{2} \sum_{i=1}^m\sum_{\ell=1}^k\|x^\ell-a^i\|^2 + \frac{\tau q}{2} \sum_{\ell=1}^k\|x^\ell\|^2\Big)\\
		&-\Big(\frac{1}{2} \sum_{i=1}^m \max_{r=1, \ldots,k}\sum_{\ell=1, \ell\neq r}^k (\|x^\ell-a^i\|)^2 +\frac{\tau}{2} \sum_{\ell=1}^k \sum_{i=1}^q \varphi_{{\Omega}_i^{\ell}}(x^\ell) \Big).
	\end{align*}
	Now, let $f=g-h$ by denoting
	\begin{align*}
		&g_1(x^1, \ldots, x^k)=\frac{1}{2} \sum_{i=1}^{m}\sum_{\ell =1}^{k}\|x^\ell - a^i\|^2, &&g_2(x^1, \ldots, x^k)= \frac{\tau q}{2}\sum_{\ell =1}^{k}\|x^\ell \|^2~,\\
		& h_1(x^1, \ldots, x^k) = \frac{1}{2} \sum_{i=1}^m \max_{r=1, \ldots,k}\sum_{\ell=1, \ell\neq r}^k \|x^\ell-a^i\|^2, &&h_2(x^1, \ldots, x^k) =\frac{\tau}{2} \sum_{\ell=1}^k \sum_{i=1}^q \varphi_{{\Omega}_i^{\ell}}(x^\ell),\end{align*}
	and let $g=g_1+g_2$ and $h=h_1+h_2$. As discussed earlier, we shall collect $x^j$ into the variable matrix $\bX$, $a^i$ into the data matrix $\bA$, and let $\bO_i={\O_i^1}\times{\O_i^2}\times\ldots\times{\O_i^k} \in \R^{k\times d}$ for $i=1,\ldots,q$. 
	
	It is clear that $g$  is differentiable and its gradient is given by
	\begin{align*}
		\nabla g(\bX) &=\nabla g_{1}(\bX)+\nabla g_{2}(\bX) =(m+\tau q)\bX-\bE \bA.
	\end{align*}
	Here, $\bE\in \R^{k\times m}$ is the matrix of ones. More details can be found in \cite{Nam2018}. Using \cref{fcg}, we find 
	\begin{equation*}
		\bX =\frac{ \mathbf Y+\bE\bA}{m+\tau q}~\in \partial g^*(\mathbf Y).
	\end{equation*}
	Next, we compute $\mathbf Y_p\in \partial h(\bX_p)$ and obtain $\bX_{p+1}$. In \cite{Nam2018}, $\mathbf W\in\partial h_1(\bX)$ is given by
	\begin{align} \label{grad2}
		\mathbf W=\sum_{i=1}^m\Big(\bX-\bA_i-e_{r(i)}(x^{r(i)}-a^i)\Big)=m\bX-\bE\bA-\sum_{i=1}^me_{r(i)}(x^{r(i)}-a^i),
	\end{align}
	where $\bA_i\in \R^{k\times d}$ is the matrix whose all rows are $a^i$, $r(i)$ is an index where the max happens for each $i$, and $e_r$ is the $k\times1$ column vector with a one in the $r^{th}$ position and zeros otherwise. For $h_2$, we have  $\mathbf U=\frac1\tau\nabla h_{2}(\bX)$ is the $k\times d$ matrix whose rows are
	$u^j=\sum_{i=1}^qP(x^j;\Omega_i^j)$. Now, let $\mathbf Y_p=\mathbf W+\tau \mathbf U$, we are able to get $\mathbf Y_p\in\partial h(\bX_p)$ at the $p^{th}$ iteration. Hence, an explicit formula for $\bX$ is 
	\begin{align*}
		\bX_{p+1}=\frac{1}{m+\tau q}\Big(m\bX_p+\tau \mathbf U-\sum_{i=1}^me_{r(i)}\big(x_p^{r(i)}-a^i\big)\Big),
	\end{align*}
	where $x_p^\ell$ denotes the $\ell^{th}$ row of $\bX_p$.
	
The summary of the DCA-based procedure can be found in  \cite[Algorithm 2]{Nam2018}. We also discussed the sensitivity of $\tau$ and provided the multiple scalar $\sigma$ and the max value $\tau_f$.  In \cite[Algorithm 3]{Nam2018} , an adaptive version of DCA when combined with a penalty parameter i s demonstrated.
	
	In \cite{Nam2018} $g$ was shown differentiable , and the subgradient of $h$ can be explicitly calculated, hence fulfilling the assumptions of the BDCA. Therefore, a BDCA-based procedure for solving \cref{PenPointclustering} is shown in \cref{BDCA2}.
	
	\begin{algorithm}[h]
		\caption{BDCA for solving \cref{Pointclustering}} \label{BDCA2}
		\begin{algorithmic}
			\Procedure{BDCA for \cref{Pointclustering} \;}{$\mathbf{A}, \mathbf{X}_0 , \Omega_{j, l}, tol, \tau, \sigma, \tau_f, \alpha > 0, \beta \in (0, 1) \,, p=1$}
			\While{$\tau < \tau_f$}
			\While{tol is FALSE}
			\State \underline{Step 1:} Find $y_p$ by executing the following iteration of DCA
			\For{$i=1,\dots,m$}
			\State Find $r(i)$ s.t. $ \|x_{p-1}^{r(i)} -a^i\|^2=\min\{\norm{x^\ell_{p-1}-a^i}^2 |\; \ell=1,\ldots,k\}$
			\State Set $\mathbf W_i:=e_{r(i)}(x_{p-1}^{r(i)}-a^i)$
			\EndFor
			\For{$\ell=1,\dots,k$}
			\State Find $u^\ell:=\sum_{j=1}^q P(x^\ell_{p-1};\O_j^\ell)$
			\EndFor
			\State Set $y_{p}:=\frac{1}{m+\tau q}\Big(m\bX_{p-1}+\tau \mathbf U-\sum_{i=1}^m \mathbf W_i\Big)$
			\State \underline{Step 2:} Set $d_p := y_p - \bX_{p-1}$
			\If{$d_p = 0$}
			\State \Return $\bX_p$
			\Else
			\State Go to \underline{Step 3}
			\EndIf
			\State \underline{Step 3:} Choose any $\Bar{\lambda}_p \geq 0$, set $\lambda_p =\Bar{\lambda}_p$
			\While{$f(y_p + \lambda_p d_p) > f (y_p) - \alpha \lambda_p^2 \Vert d_p \Vert^2$}
			\State $\lambda_p = \beta \lambda_p$
			\EndWhile
			\State  Set  $\bX_{p} = y_p + \lambda_p d_p$, and $p = p + 1$.
			\EndWhile
			\State Reassign $\tau := \sigma \tau$.
			\EndWhile
			\EndProcedure
			\State \Output{$\bX_p$}
		\end{algorithmic}
	\end{algorithm}
	
	\section{Set Clustering with Constraints}\label{Chapter:SetClusters}
	
	In this section, we revisit a model of \emph{set clustering} with constraints in \cite{Nam2018}. Given $m$ subsets $\Lambda_1,\ldots,\Lambda_m\subset\R^d$, we look for $k$ cluster centers $x^\ell\in\bigcap_{j=1}^q\O_j^\ell$ for $\ell=1,\ldots, k$, where each  $\O_j^\ell$ is a subset of $\R^d$. We also use the same square distance functions to measure distances to the sets involved. The problem of interest is given by
	\begin{equation}\label{Setclustering}
		\begin{array}{ll}
			\mbox{\rm min } &\psi(x^1, \ldots, x^k)=\sum_{i=1}^m \min_{\ell=1,\dots,k}[d(x^\ell;\Lambda_i)]^2\\
			\mbox{\rm subject to } &x^\ell\in\bigcap_{j=1}^q \O_j^\ell\; \mbox{\rm for }\ell=1, \ldots, k.\end{array}
	\end{equation}
	
	Here we assume that $\Lambda_i$ for $i=1, \ldots, m$ and $\O_j^\ell$ for $j=1,\ldots,q$ and $\ell=1, \ldots, k$ are nonempty, closed and convex.
	
	Applying the penalty method with a parameter $\tau>0$, we obtain the unconstrained set clustering problem
	\begin{equation}\label{PenSetclustering}
		\begin{array}{ll}
			\mbox{\rm min } &f(x^1, \ldots, x^k)=\frac12\sum_{i=1}^m \min\limits_{\ell=1,\dots,k}[d(x^\ell;\Lambda_i)]^2+\frac\tau2\sum_{\ell=1}^k\sum_{j=1}^q [d(x^\ell;\Omega_j^\ell)]^2,\\
			&x^1, \ldots, x^k\in \R^d.\end{array}
	\end{equation}
	Similar to the previous section, a DC decomposition of $f=g-h$ is achieved using the \emph{minimum-sum principle} and \cref{dcr} as follows
	\begin{align*}
		&g_1(\bX)=\frac{m}{2}\|\bX\|_F^2, &&g_2(\bX)=\frac{\tau q}{2}\norm \bX_F^2,\\
		&h_1(\bX)=\sum_{i=1}^m\Big (\frac{1}{2}\sum_{\ell=1}^k\varphi_{\Lambda_i}(x^\ell)+\frac{1}{2}\max_{r=1,\dots,k}\sum_{\ell=1, \ell\neq r}^k [d(x^\ell;\Lambda_i)]^2\Big),&&h_2(\bX)=\frac\tau2\sum_{\ell=1}^k \sum_{j=1}^q\varphi_{\Omega_j^\ell}(x^\ell),
	\end{align*}
	where $g=g_1+g_2$ and $h=h_1+h_2$ are convex.
	
	More details about the derivation of $f$ and computation that will be introduced later can be found in \cite{Nam2018}. Firstly, using \cref{fcg}, we can easily compute $\bX=\frac{1}{m+\tau q}\mathbf Y  \in \partial g^*(\mathbf Y)$. Then, we find $\mathbf Y\in \partial h(\bX)$ as $\mathbf Y=\mathbf V+\mathbf U$, where $\mathbf V\in \partial h_1(\bX)$ and $\mathbf U\in \partial h_2(\bX)$. Here, $h_2$ is differentiable and $\nabla h_2(\bX)=\tau \mathbf U$, where  $\mathbf U$ is the $k\times d$ matrix whose $\ell^{th}$ row is $\sum_{j=1}^kP(x^\ell;\Omega_j^\ell)$ for $\ell=1,\dots,k$. Notice that $h_1$ is not differentiable and the subgradient $\mathbf V\in \partial h_1(\bX)$ is given by 
	\[
	\mathbf V=m\bX-\sum_{i=1}^m e_{r(i)}\Big(x^{r(i)}-P(x^{r(i)}; \Lambda_i)\Big).
	\]
	As discussed in \cite{Nam2018}, for each $i =1, \ldots, m$, we choose an index $r(i)$ such that the max happens, and  $e_r$ is the $k\times1$ column vector with a one in the $r^{th}$ position and zeros otherwise. Hence, $\bX_{p+1}$ is represented by
	\begin{equation*}
		\bX_{p+1}=\frac{1}{\tau q +m}\Big(m\bX_{p}+\tau \mathbf U_p-\sum_{i=1}^m e_{r(i)}\big(x_{p}^{r(i)}-P(x_{p}^{r(i)}; \Lambda_i)\big)\Big).
	\end{equation*}
	
	In \cite{Nam2018}, the DCA-based algorithm and the adaptive $\tau$ version for solving \cref{Setclustering} were discussed in Algorithm 4 and 5. Now, we introduce the BDCA version as follows: 
	
	\begin{algorithm}[h]
		\caption{BDCA for solving \cref{Setclustering}} \label{BDCA3}
		\begin{algorithmic}
			\Procedure{BDCA for \cref{Setclustering} \; }{$\bX_0, \Lambda_i, \{\O_j^\ell\}_{j=1,\ldots,q}^{\ell=1,\ldots,k}, tol, \tau, \sigma, \tau_f, \alpha > 0, \beta \in (0, 1) \,, p=1$}
			\While{$\tau < \tau_f$}
			\While{tol is FALSE}
			\State \underline{Step 1:} Find $y_p$ by executing the following iteration of DCA
			\For{$i=1,\dots,m$}
			\For{$\ell=1,\dots,k$}
			\State Set $w_i^\ell:=P(x_{p-1}^\ell;\Lambda_i)$
			\EndFor
			\State Find $r(i)$ s.t.	$\|x^{r(i)}_{p-1}-w_i^{r(i)}\|^2=\min \limits_{\ell=1, \ldots, k}\|x^\ell_{p-1}-w_i^\ell\|^2$
			\EndFor
			\For{$\ell=1,\dots,k$}
			\State Find $u^\ell:=\sum_{j=1}^q P(x^\ell_{p-1};\O_j^\ell)$
			\EndFor
			\State Set $y_{p}:=\frac{1}{\tau q +m}\Big (m\bX_{p-1}+\tau \mathbf U_p-\sum_{i=1}^m e_{r(i)}\big(x_{p-1}^{r(i)}-w_i^{r(i)}\big)\Big)$
			\State \underline{Step 2:} Set $d_p := y_p - \bX_{p-1}$
			\If{$d_p = 0$}
			\State \Return $\bX_p$
			\Else
			\State Go to \underline{Step 3}
			\EndIf
			\State \underline{Step 3:} Choose any $\Bar{\lambda}_p \geq 0$, set $\lambda_p =\Bar{\lambda}_p$
			\While{$f(y_p + \lambda_p d_p) > f (y_p) - \alpha \lambda_p^2 \Vert d_p \Vert^2$}
			\State $\lambda_p = \beta \lambda_p$
			\EndWhile
			\State Set  $\bX_{p} = y_p + \lambda_p d_p$, and $p = p + 1$
			\EndWhile
			\State Reassign $\tau := \sigma \tau$
			\EndWhile
			\EndProcedure
			\State \Output{$\bX_p$}
		\end{algorithmic}
	\end{algorithm}
	
	\section{Numerical Experiments}\label{examples}
	
  We now implement and test the proposed algorithms in a number of examples. All the test are implemented in MATLAB R2023b, and we made use of profiling to create efficient, more realistic implementations. The code for the examples, along with the data used in generating the figures for the examples, can be found in the public GitHub \href{https://github.com/TuyentdTran/BDCAClustering.git}{github.com/TuyentdTran/BDCAClustering.git}. We run our tests on a iMac with 3.8 Ghz 8-Core i7 processor and 32 GB DDR4 memory. Throughout the examples, \cref{BDCA2,BDCA3} are tested with $\alpha = 0.05,\,\beta = 0.1,\,\tau=1,\,\sigma=10,\,\tau_f=10^8$ and the tolerance for the DCA step is $10^{-6}$, which means we terminate that step whenever $\|\bX_{p+1}-\bX_p\|_F<10^{-6}$.
	
	In examples 2 and onward, due to the size and the complexity of the problems, it is beneficial to  use the following strategy for choosing the trial step size in Step~3 of BDCA, which utilizes the previous step sizes. More details can be found in \cite{Artacho2020}.  Recall that that in \cite{AragonArtacho2018}, $\overline{\lambda}_k$ was chosen constantly equal to some fixed parameter $\overline{\lambda}>0$, we use the same strategy, choosing \(\overline{\lambda}_{p} = 2\) in \cref{BDCA2,BDCA3} for all examples except for \cref{ex76}. \medskip
	\begin{center}
		\bgroup 	\renewcommand\theenumi{\arabic{enumi}.} 	 \renewcommand\labelenumi{\theenumi}
		\fbox{%
			\begin{minipage}[t]{.98\textwidth}%
				\textbf{Self-adaptive trial step size}\medskip\\
				Fix $\gamma>1$. Set $\overline{\lambda}_0=0$. Choose some $\overline{\lambda}_1>0$ and obtain $\lambda_1$ by Step~3 of BDCA. \\ For any $k\geq 2$:
				\begin{enumerate}
					\item IF $\lambda_{k-2}=\overline{\lambda}_{k-2}$ AND $\lambda_{k-1}=\overline{\lambda}_{k-1}$ THEN set $\overline{\lambda}_{k}:=\gamma\lambda_{k-1}$;
					ELSE set $\overline{\lambda}_{k}:=\lambda_{k-1}$.
					\item Obtain $\lambda_{k}$ from $\overline{\lambda}_{k}$ by Step~4 of BDCA.
					\medskip{}
				\end{enumerate}
		\end{minipage}}
		\egroup{}
	\end{center}\medskip
	
	The \emph{self-adaptive strategy} here uses the step size that was chosen in the previous iteration as a new trial step size for the next iteration, except in the case where two consecutive trial step sizes were successful. In that case, the trial step size is increased by multiplying the previously accepted step size by $\gamma>1$.  In all our experiments we took $\gamma:=2$. Furthermore we choose  the initial step size \(\overline{\lambda}_{1}=2\), for all examples, the same as the non-adaptive BDCA. 
	
	\subsection{Constrained Clustering}\label{ssec: Constrained Clustering}
	
	\begin{Example}\label{ex76}
		We first consider the same example as\cite[Example 7.1~]{Nam2018} to compare BDCA against one of the original examples of \cite{Nam2018}. We use the dataset EIL76 taken from the Traveling Salesman Problem Library \cite{rei} and impose the following constraints on the solution:
		\begin{enumerate}
			\item The first center is a common point of a box whose vertices are $(40,40)$; $(40,60)$; $(20,60)$; $(20,40)$ and a ball of radius $r=7$ centered at $(20,60)$.
			\item The second center is in the intersection of two balls of the same radius $r=7$, centered at $(35,20)$ and $(45,22)$, respectively.
		\end{enumerate}
		The initial centers are chosen as follows:
		\begin{itemize}
			\item The first center is drawn randomly from the box.
			\item The second center is randomly chosen from the ball centered at $(35,20)$.
		\end{itemize}
		For this problem we take the trial step size \(\overline{\lambda}_{p} =1\) for BDCA. We run the test 100 times to achieve the following approximate average solutions and cost values for DCA and BDCA respectively.
		\begin{align*}
			&\text{\textbf{DCA: }}  &&\bX =
			\begin{pmatrix}
				26.69959      &57.97127\\
				41.06910	    &23.48800
			\end{pmatrix}, 
			&& \text{Cost: }\psi(\bX)=33576.25344, \\
			& \text{\textbf{BDCA: }} && \bX = 
			\begin{pmatrix}
				26.69959      &57.97125\\
				41.06910	    &23.48789
			\end{pmatrix},
			&& \text{Cost: } \psi(\bX)=33576.25387,
		\end{align*}		
		with the cost fluctuating within the range of $10^{-11}$ for both BDCA and DCA between runs.
		
		In \cref{76_compare}, we compare the ratio between time to complete DCA over BDCA and the ratio of total number of iterations for DCA over BDCA. We can see that BDCA runs about 1.5 times faster and takes about \(1/4\) of the iterations. In spite of the \(1/4\) reduction in iterations, we only see 1.5 times speed up due to the non-trivial cost of the line search. The average run time for DCA and BDCA are 0.0038s and 0.0024s respectively.
		\begin{figure}[h!]
			\centering
			\begin{subfigure}[b]{0.49\textwidth}
				\centering
				\includegraphics[width=\textwidth]{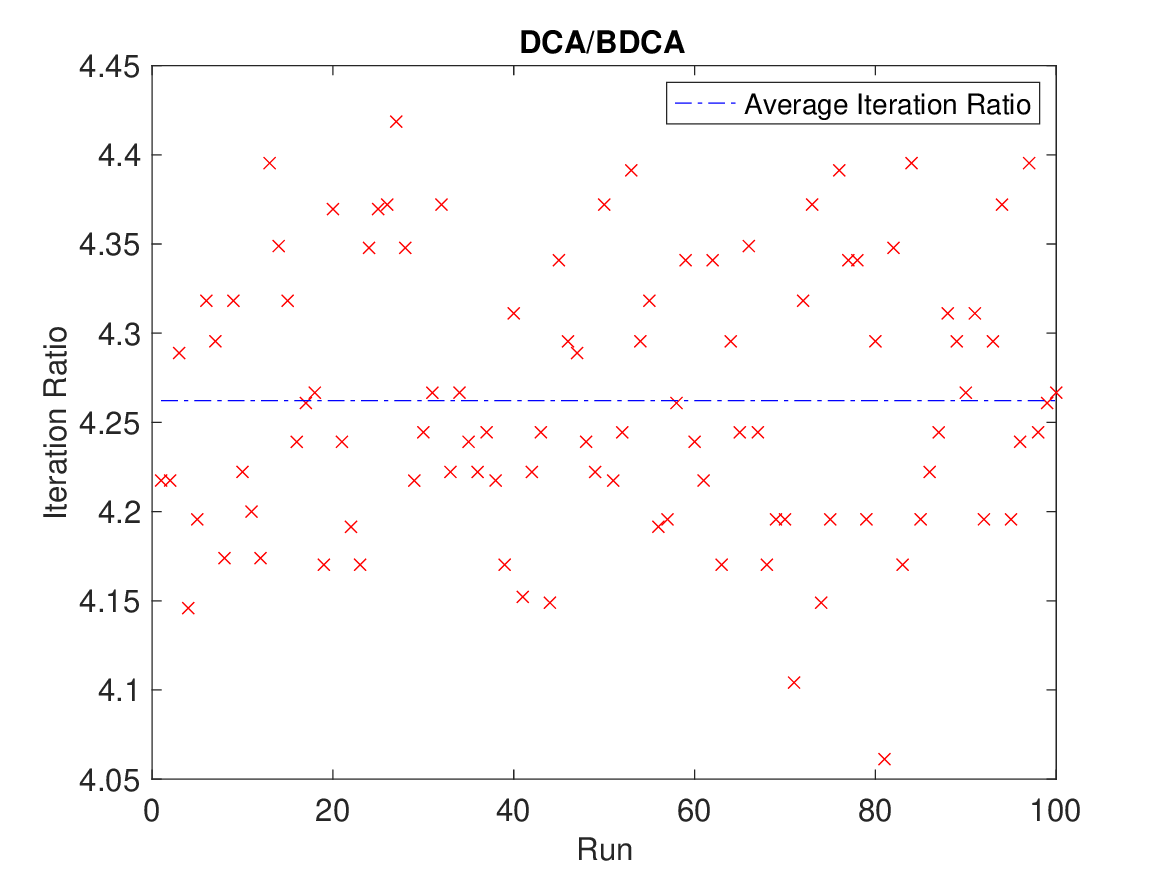}
			\end{subfigure}
			\hfill
			\begin{subfigure}[b]{0.49\textwidth}
				\centering
				\includegraphics[width=\textwidth]{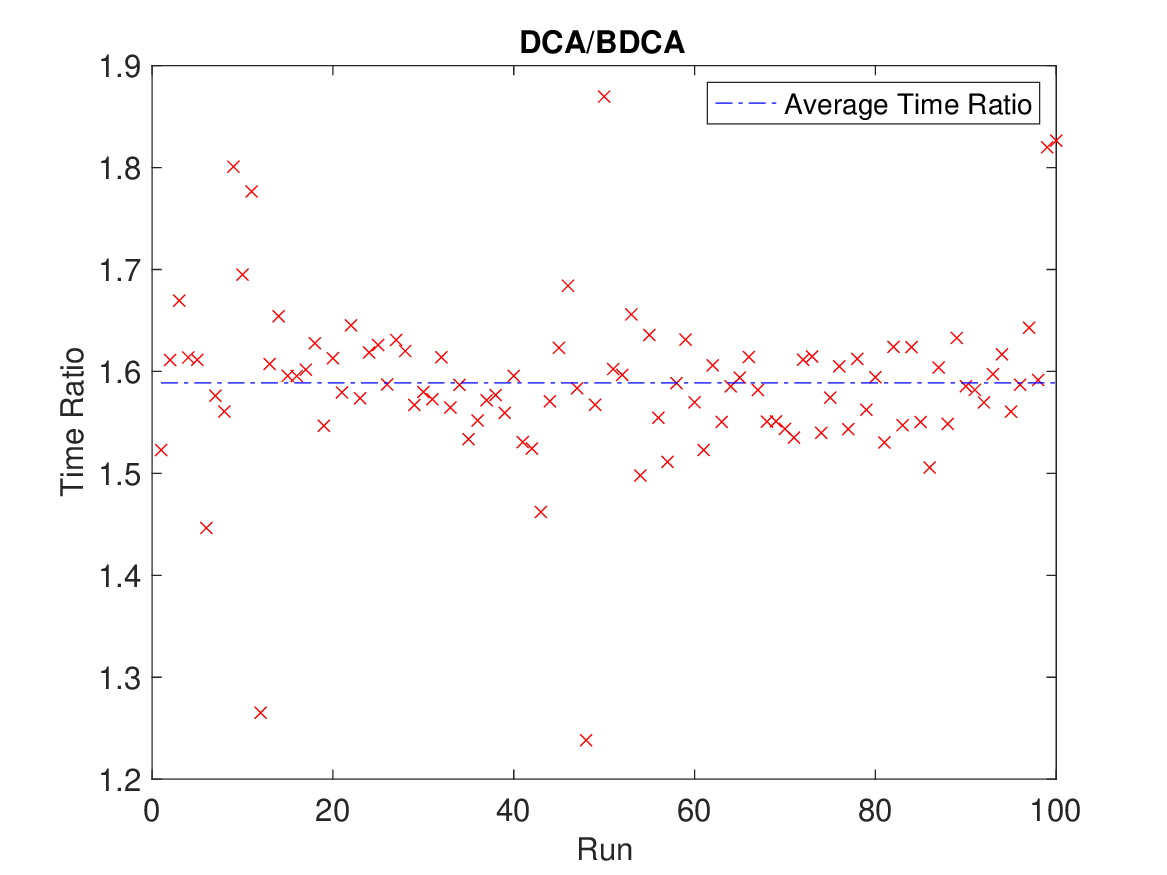}
			\end{subfigure}
			\caption{Iteration and Time Ratio Comparison for the dataset EIL76 with 2 centers, \cref{ex76}.}
			\label{76_compare}
		\end{figure}
		A visualization of the problem is demonstrated in \cref{Fig3}.
		\begin{figure}[h!]
			\centering
			\includegraphics[width=0.6\linewidth]{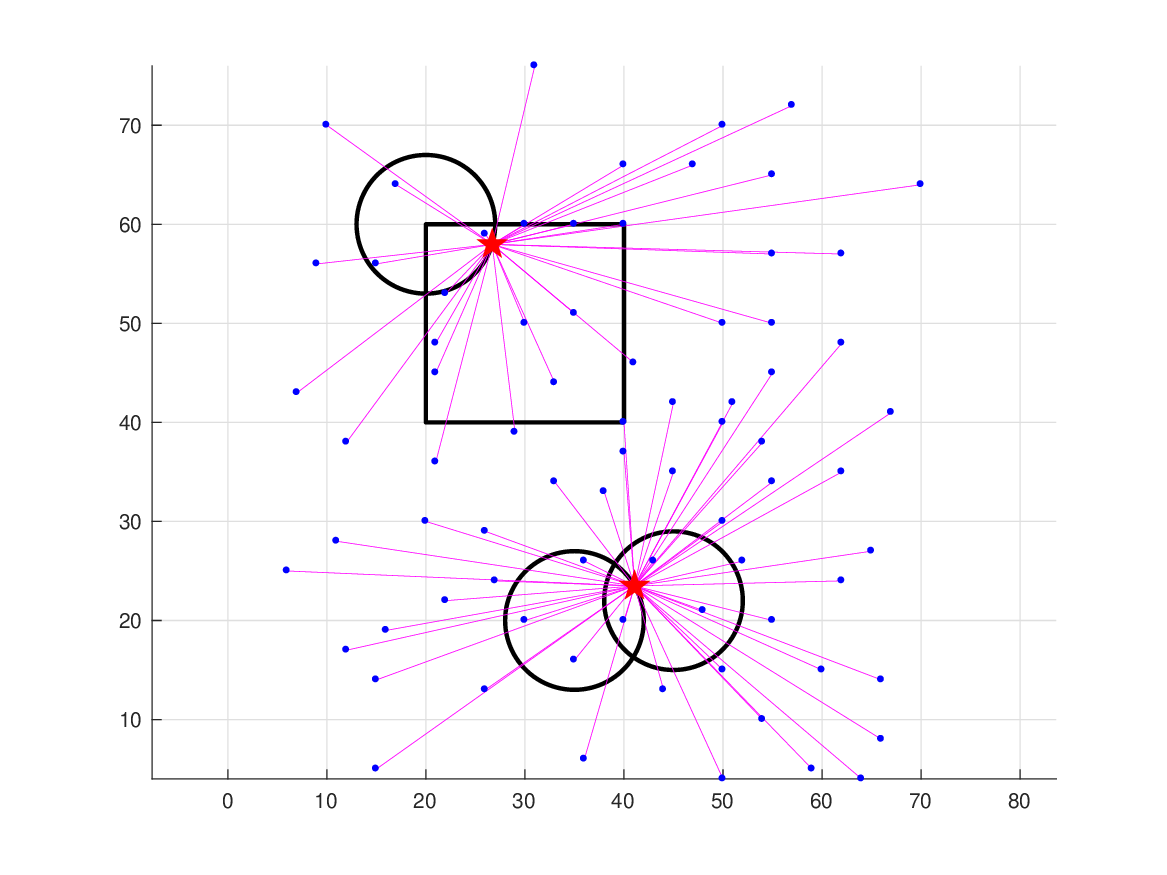}
			\caption{A 2-center constrained clustering problem for dataset EIL76, \cref{ex76}.}
			\label{Fig3}
		\end{figure}
	\end{Example}
	
	\begin{Example}\label{Example: Scaling}
		We now perform a scaling study to test the performance of BDCA vs DCA for a variety of number of points and dimensions. The intent is give a rough idea of what sort of behavior can be expected for clustering with constraints problems in different problem regimes. In this numerical experiment, we generated $n$ random points from continuous uniform distribution on the interval $[0,10]^m$. Here, $m\in\{2,3,5,10,20\}$ and $n\in\{50,100,500,1000,5000,10000,50000\}$. We impose 3 ball constraints of radius 1 with centers as follows:
		\begin{enumerate}
			\item Repeating in the pattern $[1,5,1,5\ldots]$
			\item Repeating in the pattern $[6,4,6,4,\ldots]$
			\item $[8,\ldots,8]$
		\end{enumerate}
		For each combination of $n$ and $m$,  we run with 100 random starting points drawn from the constraints and then test the performance of BDCA vs DCA.
		\begin{figure}[h!]
			\centering
			\begin{subfigure}[b]{0.49\textwidth}
				\centering
				\includegraphics[width=\textwidth]{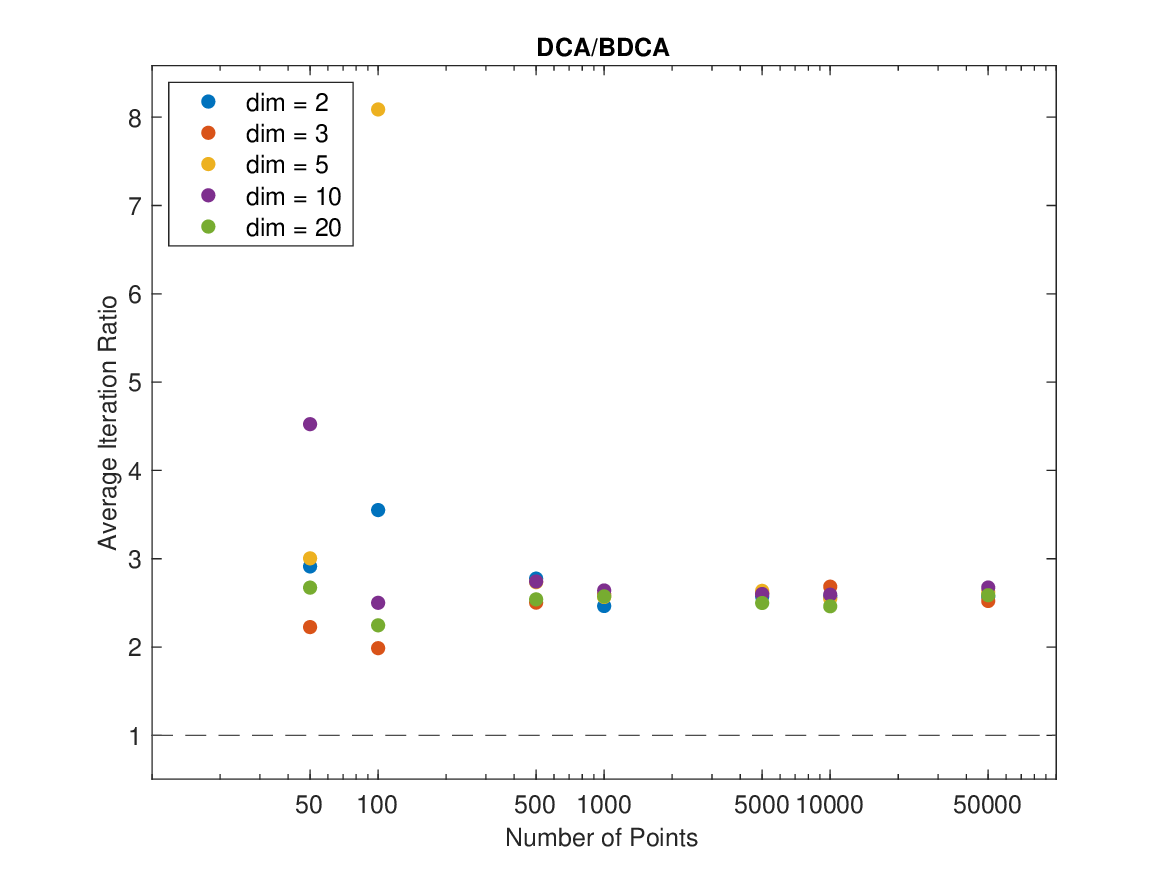}
			\end{subfigure}
			\hfill
			\begin{subfigure}[b]{0.49\textwidth}
				\centering
				\includegraphics[width=\textwidth]{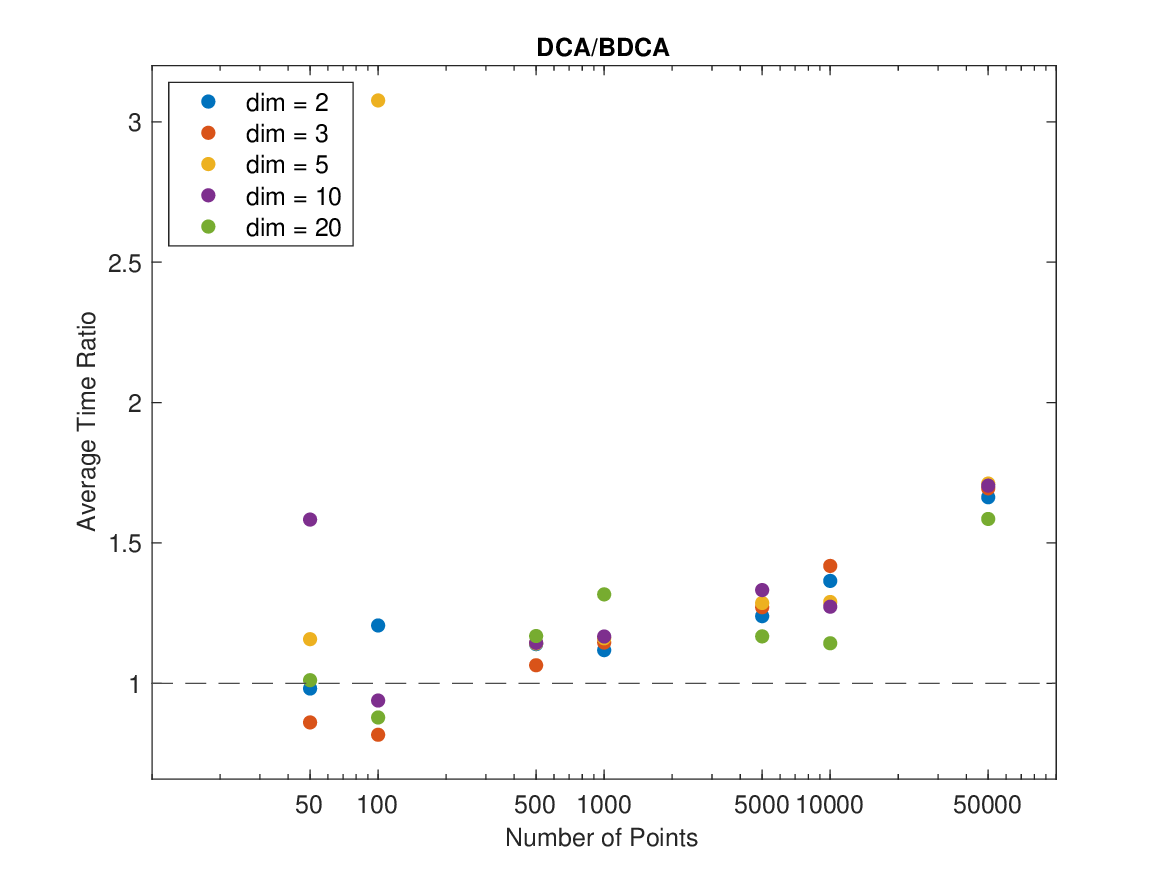}
			\end{subfigure}
			\caption{Iteration and Time Ratio Comparison between DCA and BDCA for \cref{Example: Scaling}.}
			\label{compare_BDCA}
		\end{figure}
		\begin{figure}[h!]
			\centering
			\begin{subfigure}[b]{0.49\textwidth}
				\centering
				\includegraphics[width=\textwidth]{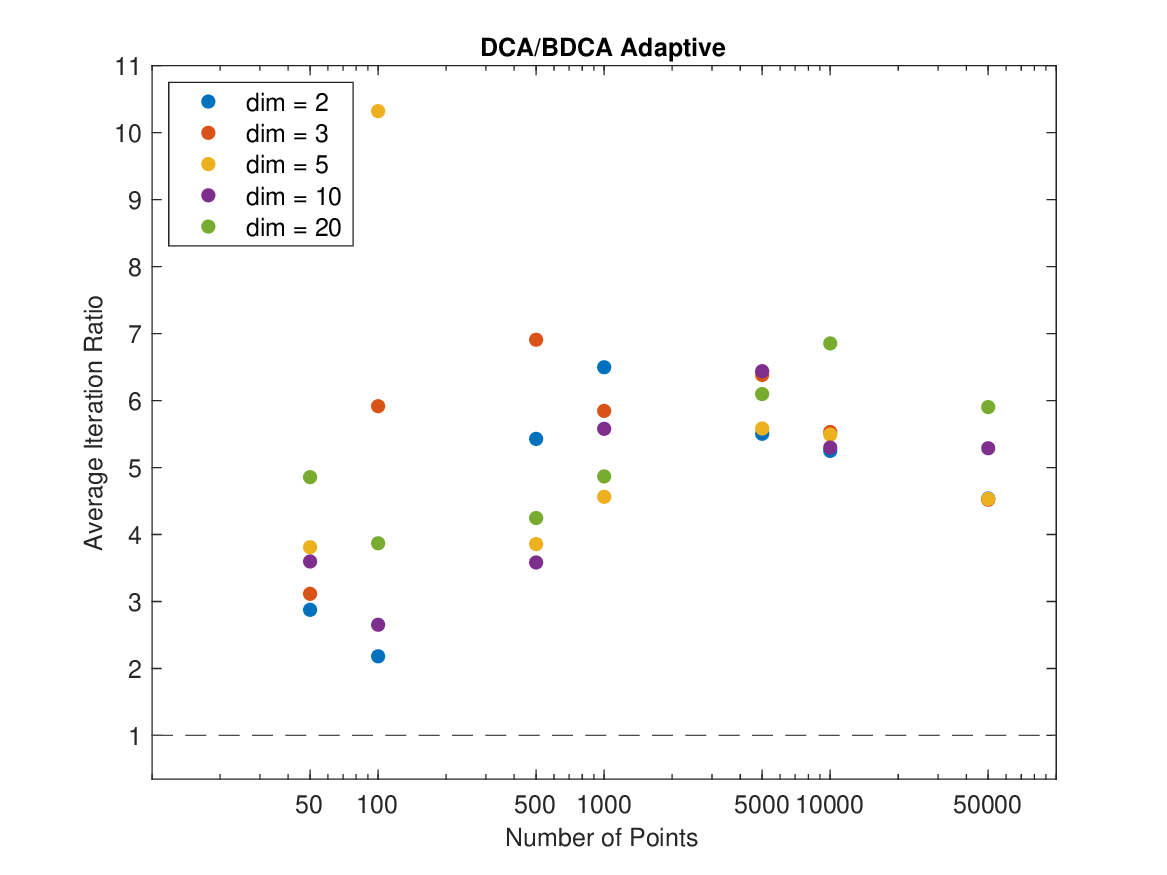}
			\end{subfigure}
			\hfill
			\begin{subfigure}[b]{0.49\textwidth}
				\centering
				\includegraphics[width=\textwidth]{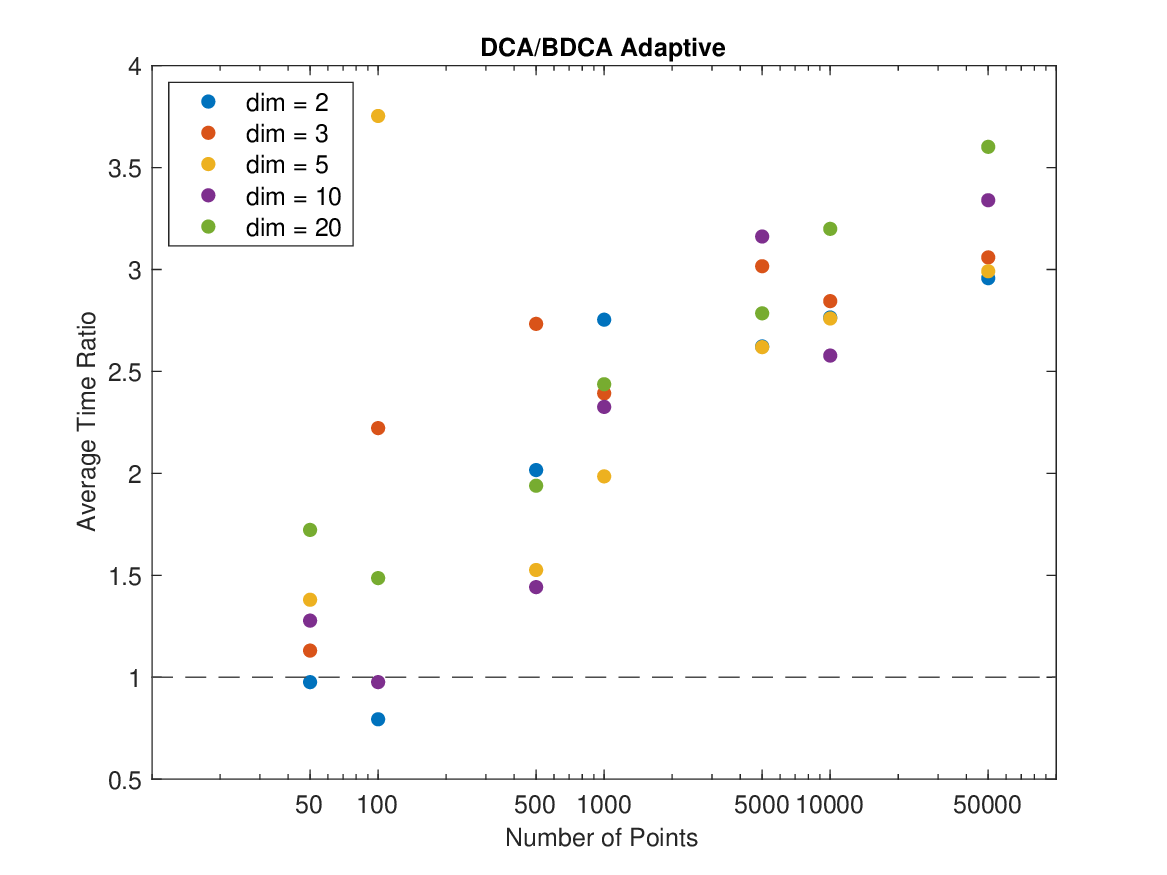}
			\end{subfigure}
			\caption{Iteration and Time Ratio Comparison between DCA and adaptive BDCA for \cref{Example: Scaling}.}
			\label{compare_BDCA_adap}
		\end{figure}
		
		\begin{figure}[h]
			\centering
			\includegraphics[width=0.7\textwidth]{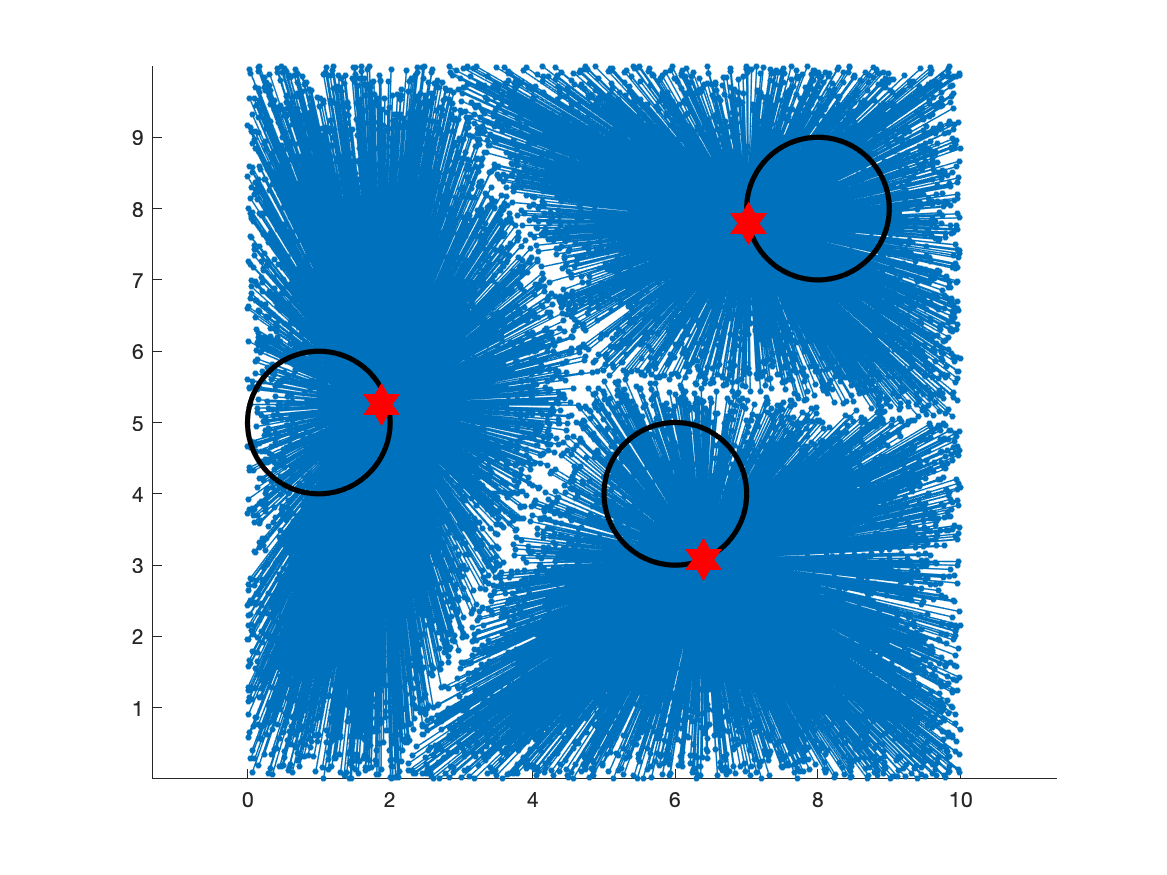}  
			\caption{A visualization of a 3-center clustering problem with 10000 points drawn from $[0,10]^2$, \cref{Example: Scaling}. Red stars are centers.}
			\label{10000pts}
		\end{figure}
		\Cref{compare_BDCA,compare_BDCA_adap} show that on average, both BDCA and adaptive BDCA are better than DCA in term of iterations and time. We can also observe that for most cases, as the dimension and number of points increase, adaptive BDCA is better than DCA for both time and iterations. On the other hand, BDCA is also better than DCA when we increase the number of points. In the case we have a small number of  points, adaptive BDCA is still always better than DCA but not the BDCA with constant value of  $\lambda$. Even though for BDCA the number of iterations is significantly less, as mentioned in \cref{ex76}, this does not come for free. You still have to compensate for the time to perform the line search, and for a small number of points, with the MATLAB implementation, the iteration reduction from BDCA isn't enough to compensate for the line search cost. The self-adaptive BDCA is necessary in this low point regime to have a time improvement. In all scenarios, iteration counts are always at least 2 times better than DCA and when the number of points go from $500$ and above, we see the improvement in time for basic BDCA. 
		
		In \cref{table: DCA,table: BDCA,table: BDCA Adaptive}, the average run time and standard deviation are reported for each situation. From those results and those of \cref{compare_BDCA,compare_BDCA_adap}, a suggested approach is to always use BDCA with self-adaptivity, except perhaps for a low number of points, where depending on your implementation and problem it may be faster to simply use DCA.
	\end{Example}
	
	\subsection{Set Clustering with Constraints}\label{ssec: set clustering}
	
	\begin{Example}\label{ex50}
		We now use \cref{BDCA3} to solve a set clustering problem with constraints which was previously discussed in \cite[Example 7.2~]{Nam2018}. We consider the latitude and longitude of the 50 most populous US cities taken from the $2014$ United States Census Bureau data \footnote{\url{https://en.wikipedia.org/wiki/List_of_United_States_cities_by_population}}, and approximate each city by a ball with radius $0.1\sqrt{\frac{A}{\pi}}$ where $A$ is the	city's reported area in square miles.
		
		We set up 3 centers as before with the requirement that each center must belong to the intersection of two balls. The centers of these constrained balls are the columns of the matrix below
		$$
		\begin{pmatrix}
			-80 &  -80 & -92& -90& -115&  -110\\
			\;\;\;34 &\;\;\;38& \;\;\;37& \;\;\;40& \;\;\;45& \;\;\;40\\
		\end{pmatrix}
		$$
		with corresponding radii given by $
		\begin{pmatrix}
			2& 3& 4& 3& 4& 4\\
		\end{pmatrix}
		$. A visualization for this problem using a plate Carr\'{e}e projection was plotted as in \cite[Example 7.2~]{Nam2018} in \cref{50US} below \footnote{\url{https://www.mathworks.com/help/map/pcarree.html}}.
		\begin{figure}[h]
			\centering
			\includegraphics[scale=0.45]{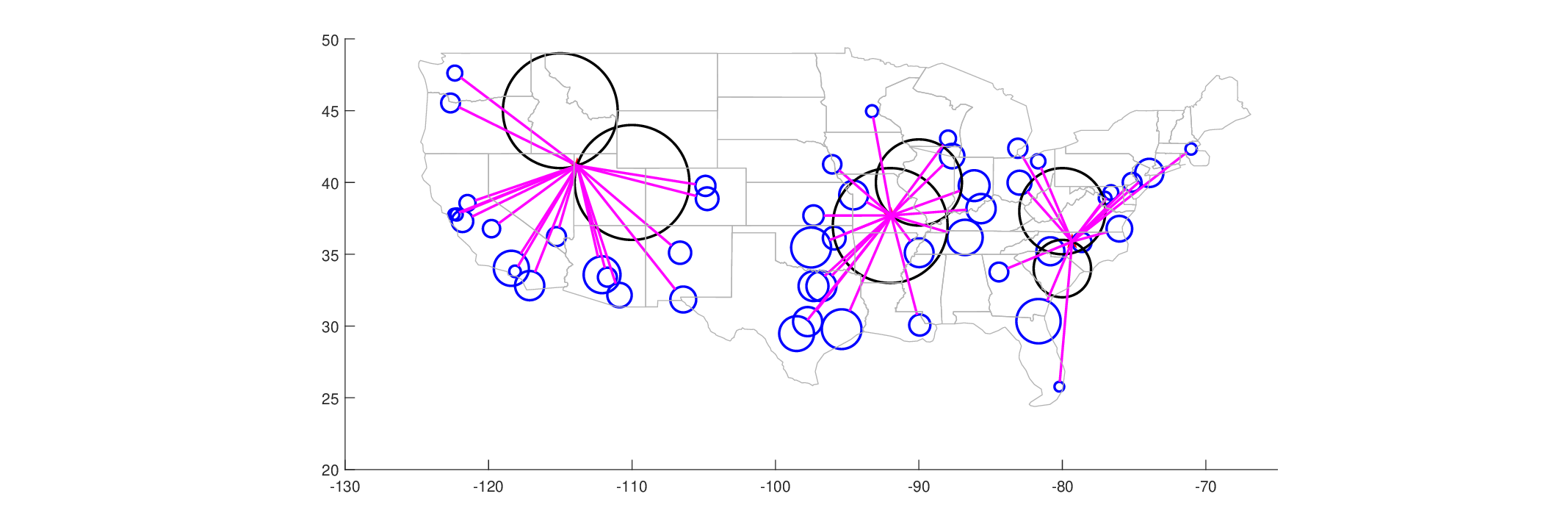}  
			\caption{A 3-center set clustering problems with 50 most populous US cities. Each city is approximated by a ball proportional to its area, \cref{ex50}.}
			\label{50US}
		\end{figure}
		In this experiment, we run the test 100 times for DCA, BDCA and adaptive BDCA. The initial centers are drawn randomly from points belonging to the first 3 constrained balls. This yields the following approximate average solutions and cost values for each run
		\begin{align*}
			&\text{\textbf{DCA: }}  &&
			\bX=\begin{pmatrix}
				-79.32232   &35.88171\\
				-91.93111   &37.70414\\
				-113.82298   &41.17699	\end{pmatrix}, 
			&& \text{Cost: }\psi(\bX)=2271.07289, \\
			& \text{\textbf{BDCA: }} && \bX = 
			\begin{pmatrix}
				-79.32301  &35.88195\\
				-91.93127   &37.70427\\
				-113.82299   &41.17698	\end{pmatrix},
			&& \text{Cost: } \psi(\bX)=2271.07019, \\
			& \text{\textbf{Adaptive BDCA: }} && \bX = 
			\begin{pmatrix}
				-79.32302  &35.88196\\
				-91.93094   &37.70400\\
				-113.82299   &41.17697	\end{pmatrix},
			&& \text{Cost: } 
			\psi(\bX)=	2271.06986.
		\end{align*}
		Note that they are equivalent up to the relative tolerance of \(10^{-6}\).
		\begin{figure}[h!]
			\centering
			\begin{subfigure}[b]{0.49\textwidth}
				\centering
				\includegraphics[width=\textwidth]{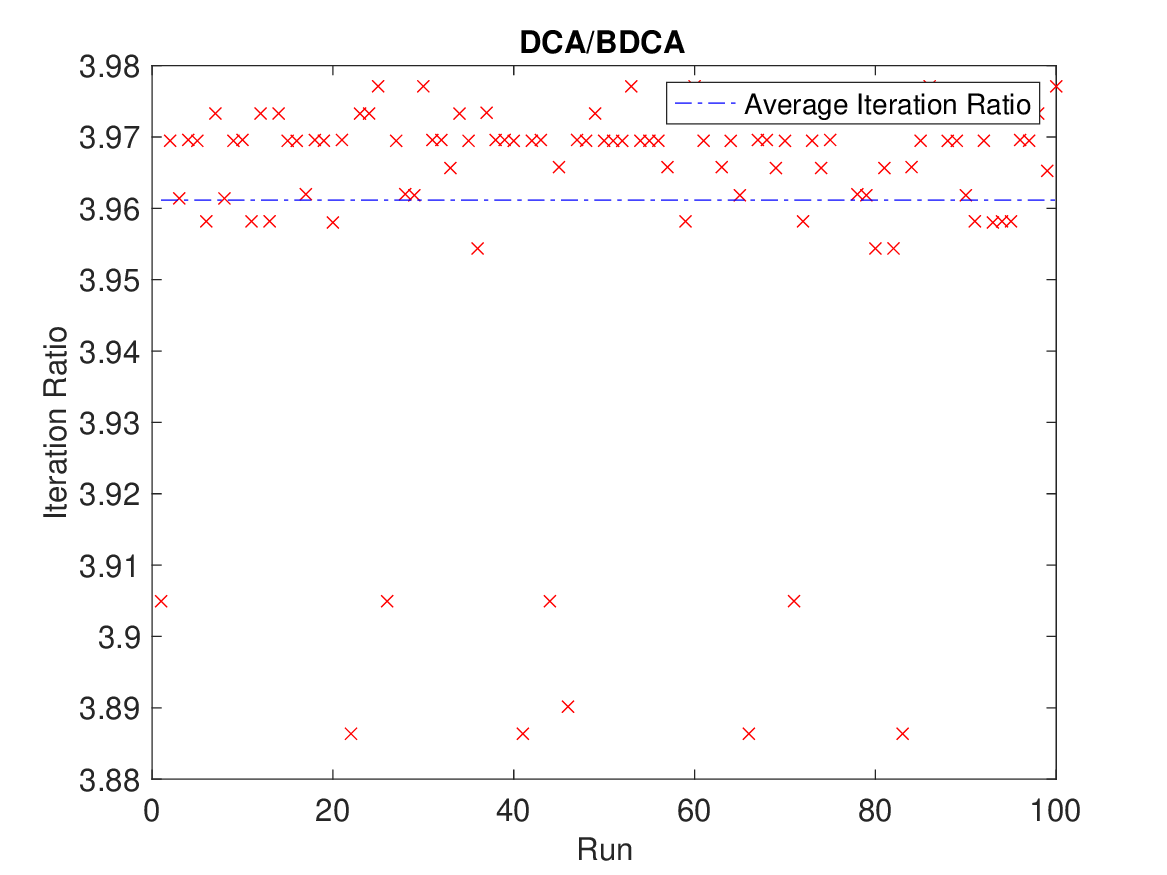}
			\end{subfigure}
			\hfill
			\begin{subfigure}[b]{0.49\textwidth}
				\centering
				\includegraphics[width=\textwidth]{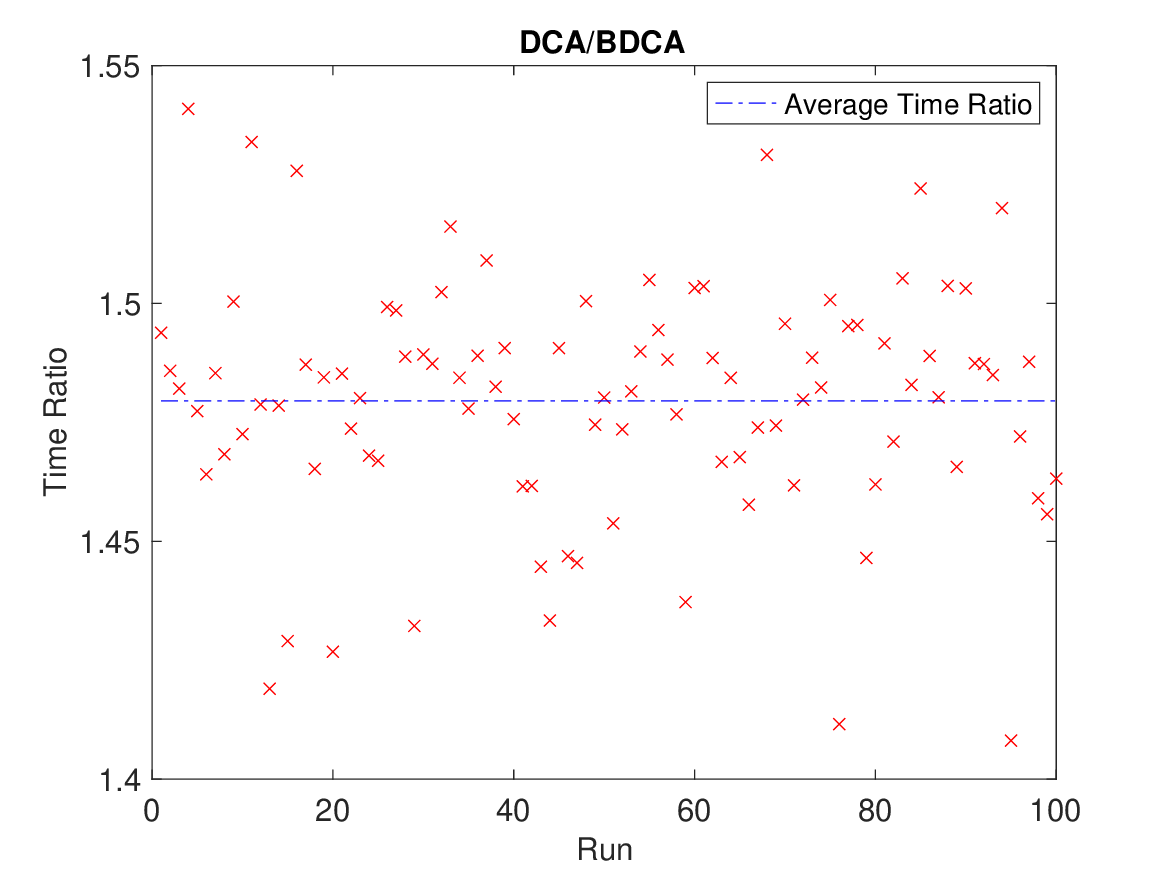}
			\end{subfigure}
			\caption{Iteration and Time Ratio Comparison between DCA and BDCA for 50 most populous US cities with 3 constraints, \cref{ex50}.}
			\label{50_compare_BDCA}
		\end{figure}
		\begin{figure}[h!]
			\centering
			\begin{subfigure}[b]{0.49\textwidth}
				\centering
				\includegraphics[width=\textwidth]{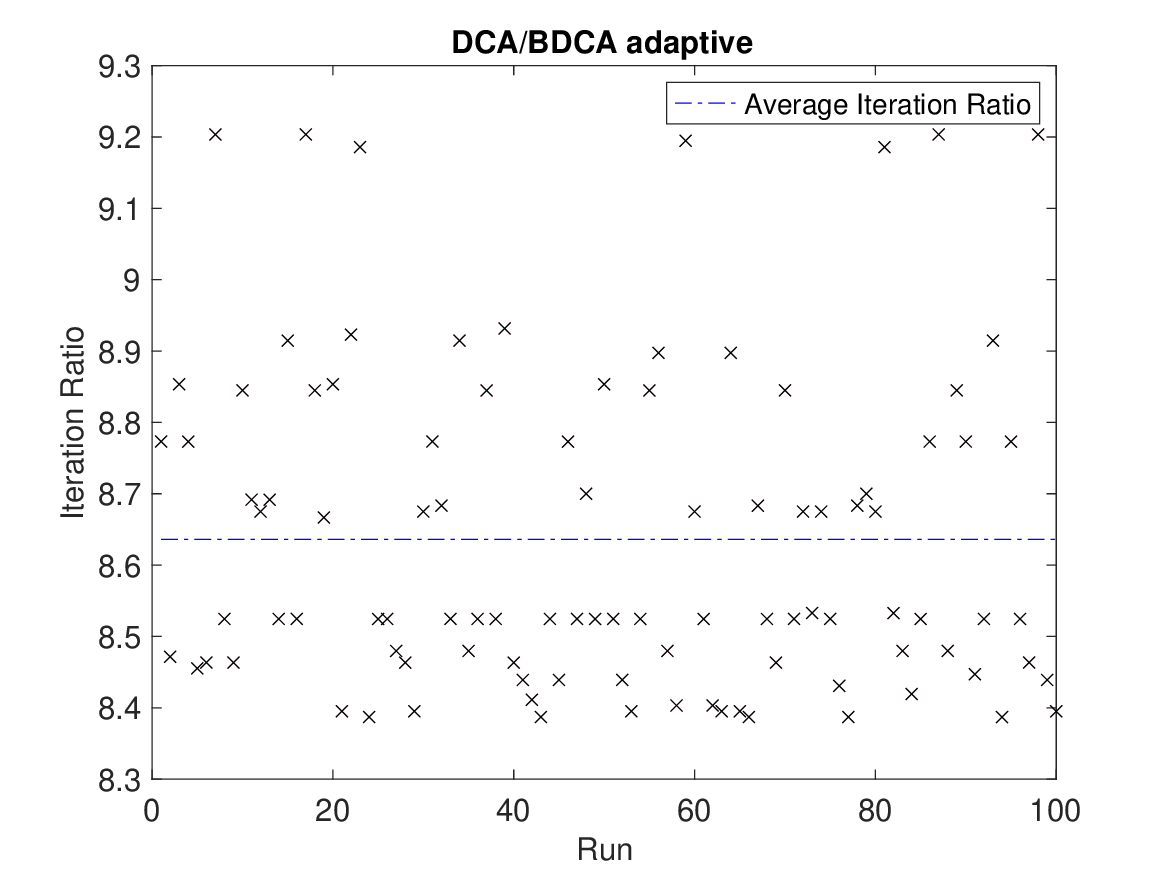}
			\end{subfigure}
			\hfill
			\begin{subfigure}[b]{0.49\textwidth}
				\centering
				\includegraphics[width=\textwidth]{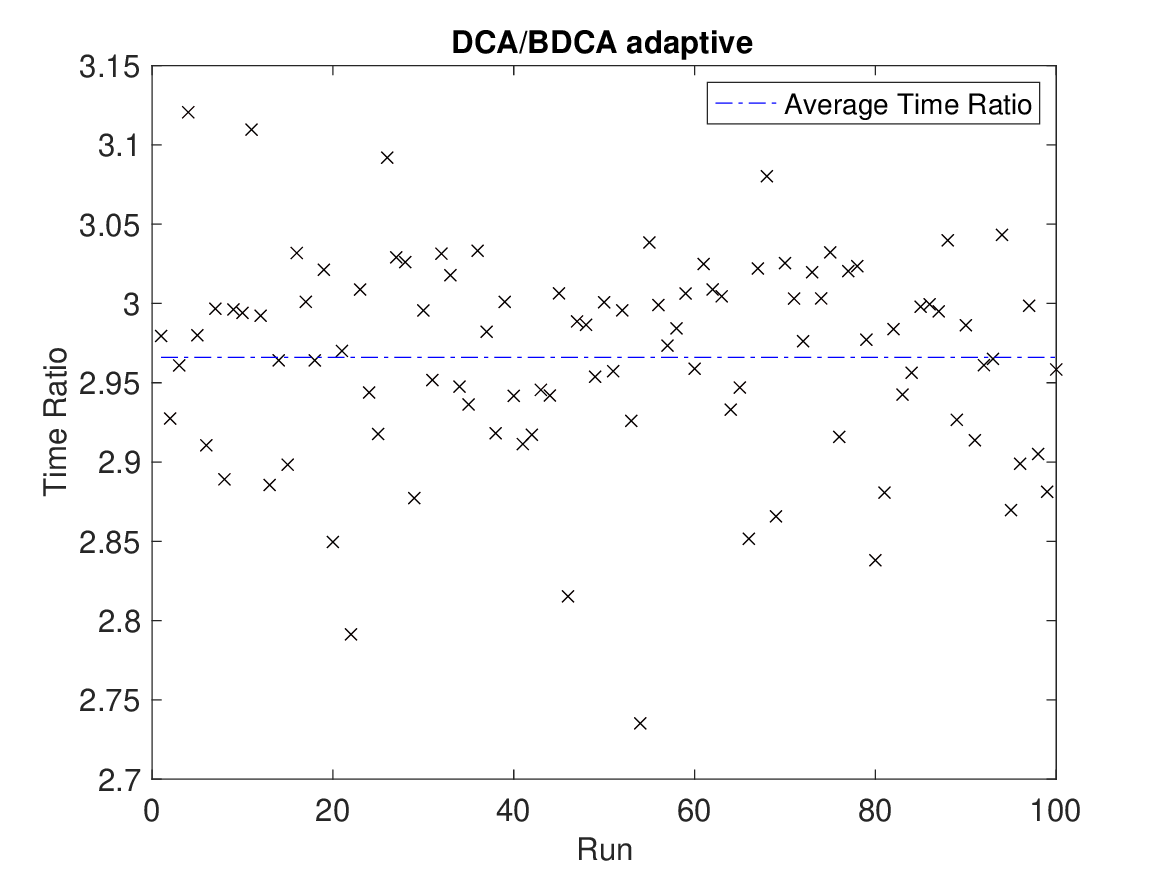}
			\end{subfigure}
			\caption{Iteration and Time Ratio Comparison between DCA and adaptive BDCA for 50 most populous US cities with 3 constraints, \cref{ex50}.}
			\label{50_compare_BDCA_adap}
		\end{figure}
		In  \cref{50_compare_BDCA}  we compare the ratio between time DCA over BDCA and the ratio of total number of iterations for DCA over BDCA. Similarly, in \cref{50_compare_BDCA_adap}, we compare the ratio between DCA over adaptive BDCA for time and number of iterations. The dashed lines in both figures show the overall average ratio for both time and iterations. We can see that adaptive BDCA outperforms BDCA and both of them are better than DCA. On average, DCA is slower than BDCA by 1.48 times and 2.95 times for adaptive BDCA. In term of iterations, DCA requires 3.96 times more than BDCA and 8.6 times for adaptive BDCA. We see here that the \emph{self-adaptive} BDCA represents a significant improvement over regular BDCA. Since the starting points are chosen randomly within constraints, we can see that the ratio for time and iterations are scattered for both comparisons.
	\end{Example}
	\begin{Example}\label{ex1500}
		We next consider the latitude and longitude of the 1500 most populous US cities derived from $2023$ United States Census Bureau data \footnote{\url{https://simplemaps.com/data/us-cities}}, and approximate each city by a ball with radius $10^{-3}\sqrt{\frac{A}{\pi}}$ where $A$ is the	city's reported area in square miles.
		We impose the following constraints on the solution: \begin{enumerate}
			\item One center is to lie within $4^\circ$ latitude/longitude of Caldwell, Idaho and inside the rectangular box with coordinates $[-115,42;-115,49;-125,49;-125,42]$.
			\item One center is to lie within the state of Colorado  and within $2.5^\circ$ latitude/longitude of Cheyenne, WY.
			\item One center is to lie  within $3^\circ$ latitude/longitude of Chicago, Illinois and  within $4^\circ$ latitude/longitude of St. Louis, MO.
			\item One center is to lie East of $-75^\circ$ longitude and within $4^\circ$ latitude/longitude of Washington,DC.
		\end{enumerate}
		This example demonstrates the ability to handle more complicated constraints than in \cref{ex50}, as well as how the algorithms scale as you consider more points when compared to \cref{ex50}.
		\begin{figure}[h!]
			\centering
			\includegraphics[width=\textwidth]{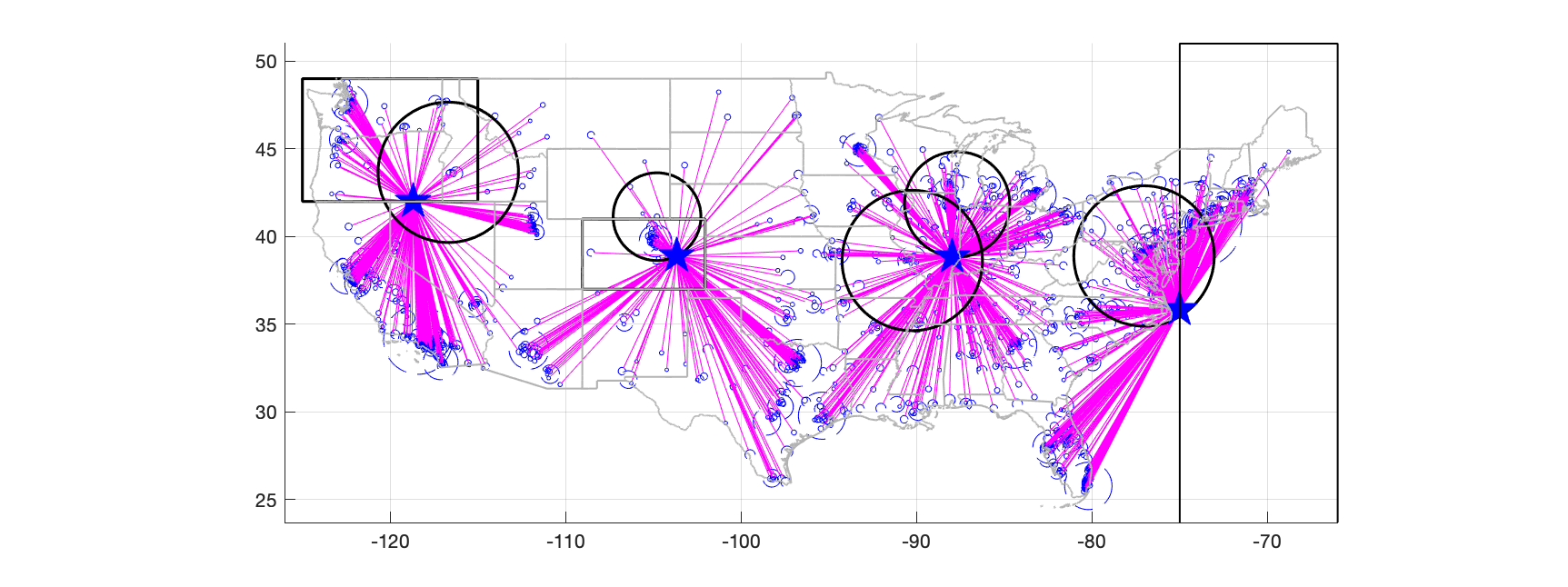}
			\caption{A 4-center set clustering with constraints for 1500 most populous  US cities, \cref{ex1500}.}
			\label{MF2}
		\end{figure}
		
		In  \cref{1500_compare_BDCA}  we compare the ratio between time DCA over BDCA and the ratio of total number of iterations for DCA over BDCA. Similarly, in \cref{1500_compare_BDCA_adap}, we compare the ratio between DCA over adaptive BDCA for time and number of iterations. The dashed lines in both figures show the overall average ratio for both time and iterations. From \cref{1500_compare_BDCA_adap,1500_compare_BDCA} we can see that BDCA improves iterations by about 5.4 times and run time by about 1.9, while adaptive BDCA gives a much better performance improvement of 12.6 for iterations and nearly 4.2 times improvement in run time. We see that compared with \cref{ex50} adaptive BDCA offers even greater improvements in runtime and iterations as the problem size increases. A trend that we will see again in the set scaling example, \cref{Example: Set Clustering}.  
	\end{Example}
	\begin{figure}[h!]
		\centering
		\begin{subfigure}[b]{0.49\textwidth}
			\centering
			\includegraphics[width=\textwidth]{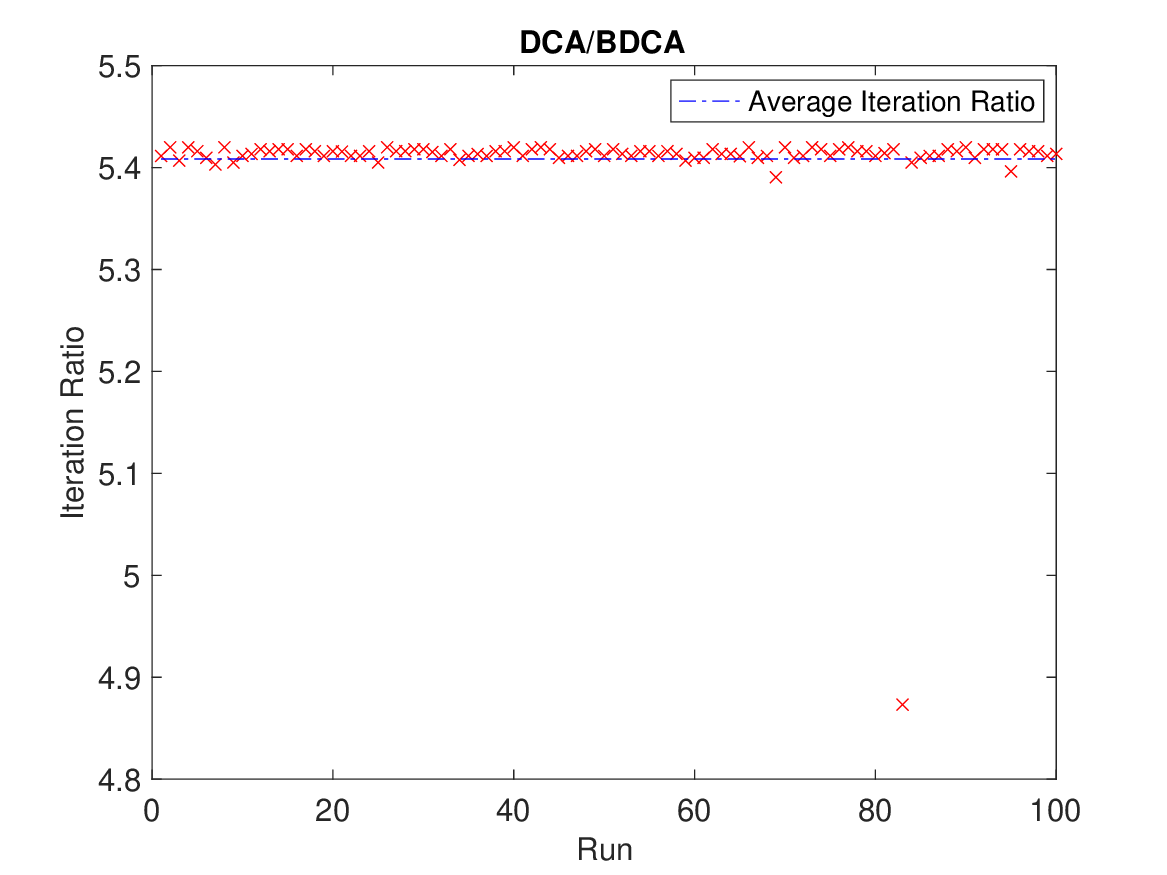}
		\end{subfigure}
		\hfill
		\begin{subfigure}[b]{0.49\textwidth}
			\centering
			\includegraphics[width=\textwidth]{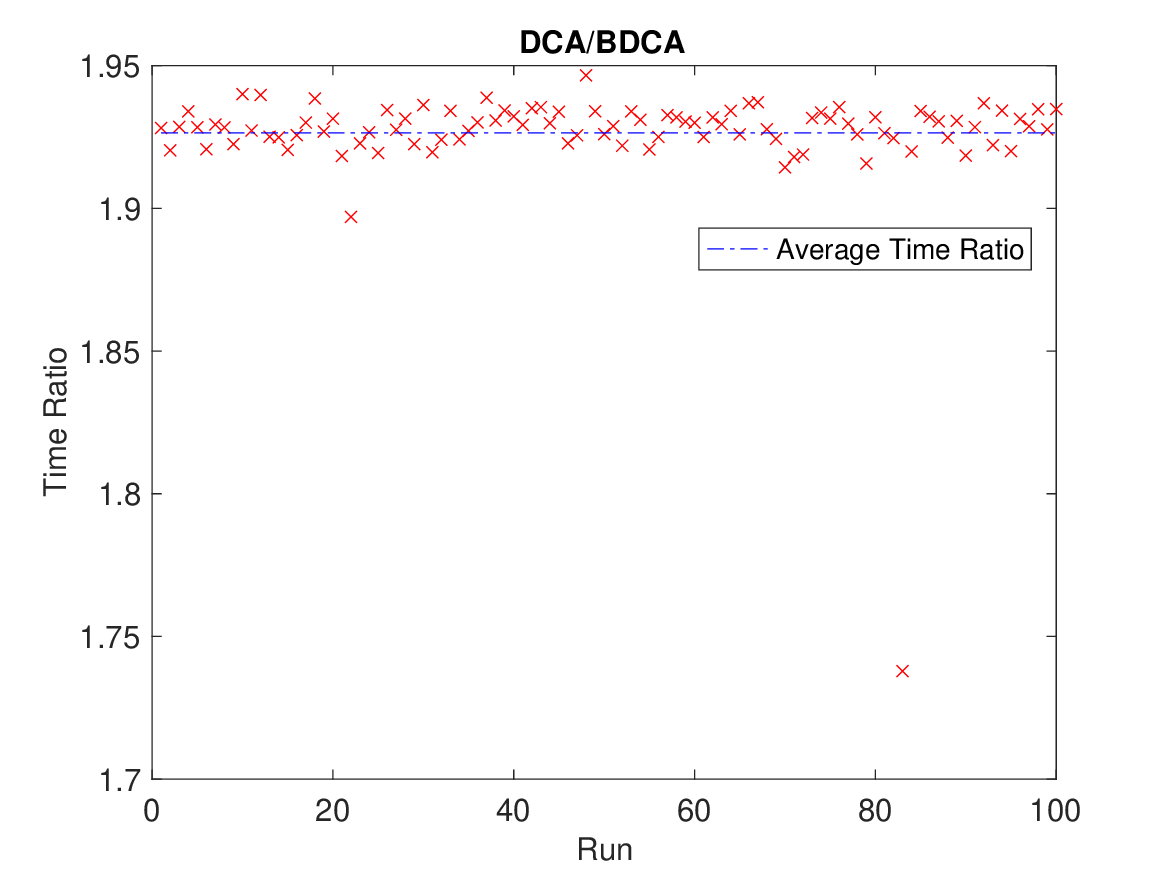}
		\end{subfigure}
		\caption{Iteration and Time Ratio Comparison between DCA and BDCA for 1500 most populous US cities with 4 constraints, \cref{ex1500}.}
		\label{1500_compare_BDCA}
	\end{figure}
	\begin{figure}[h!]
		\centering
		\begin{subfigure}[b]{0.49\textwidth}
			\centering
			\includegraphics[width=\textwidth]{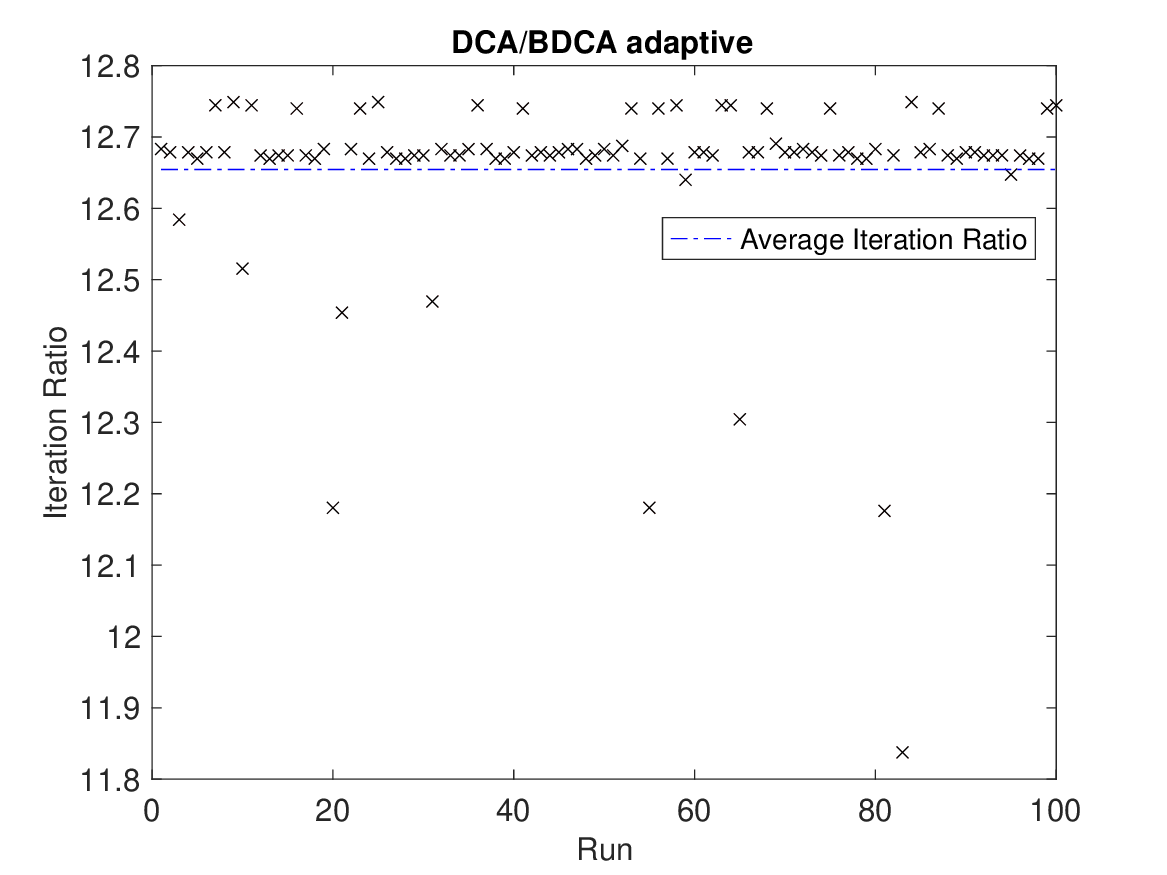}
		\end{subfigure}
		\hfill
		\begin{subfigure}[b]{0.49\textwidth}
			\centering
			\includegraphics[width=\textwidth]{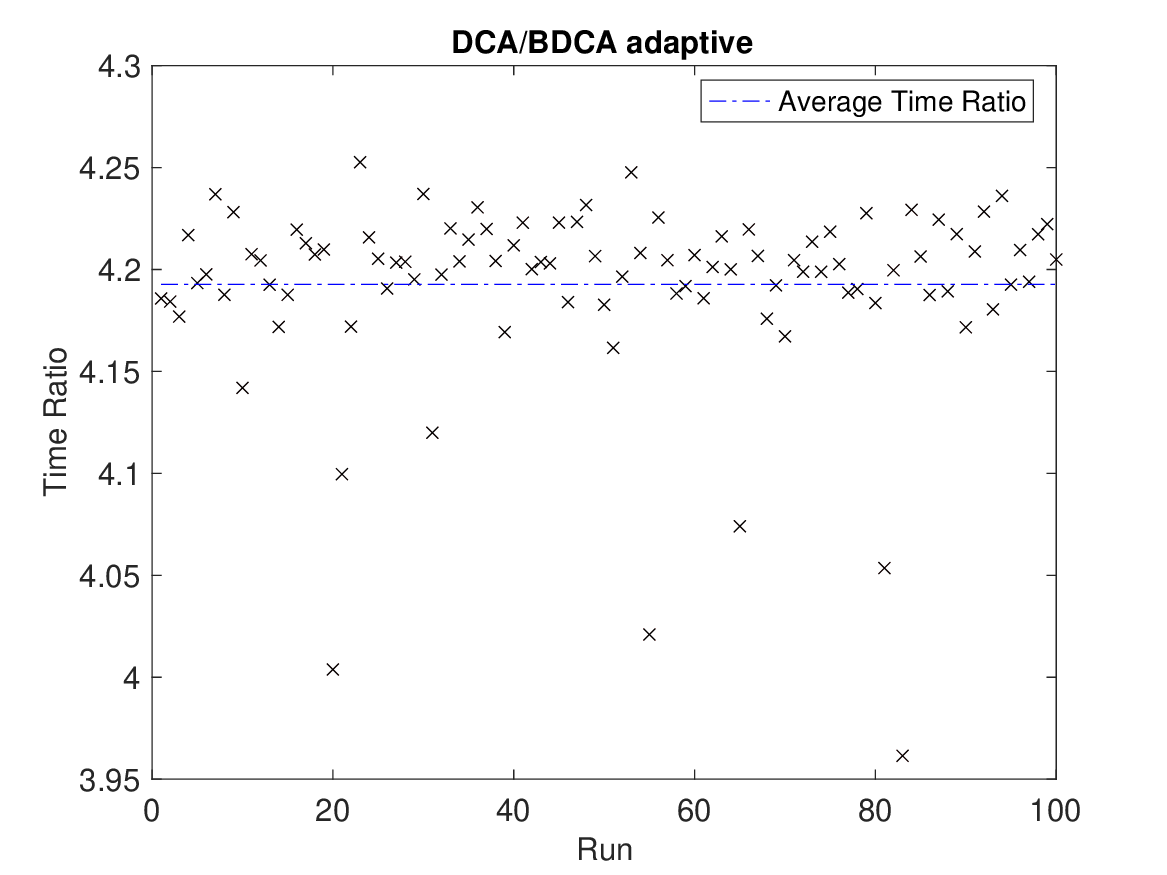}
		\end{subfigure}
		\caption{Iteration and Time Ratio Comparison between DCA and adaptive BDCA for 1500 most populous US cities with 4 constraints, \cref{ex1500}.}
		\label{1500_compare_BDCA_adap}
	\end{figure}
	\begin{Example}\label{Example: Set Clustering}
		We again perform a scaling study to test the performance of BDCA vs DCA for a variety of number of points and dimensions. The intent is give a rough idea of what sort of behavior can be expected for set clustering with constraints problems in different problem regimes. In this numerical experiment, we generated $n$ random points from continuous uniform distribution on the interval $[0,10]^m$. Here, $m\in\{2,3,5,10\}$ and $n\in\{50,100,250,500,1000,5000,10000,50000\}$. We impose 4 constraints, each formed by the intersection of two balls of radius 1 with centers composed of the first $m$ entries of the following vectors: 
		\begin{enumerate}
			\item \([1, 5,5,\hdots, 5] \text{ and }  [2,6,5,\hdots,5 ] \) 
			\item \([5, 4, 1, 2, 3, 1, 2, 3, 1, 2] \text{ and } [4, 4, 1, 2, 3, 1, 2, 3, 1, 2]\) 
			\item \([8, 5, 9, 8, 7, 9, 8, 7, 9, 8] \text{ and } [8, 4, 9, 8, 7, 9, 8, 7, 9, 8]\)  
			\item \([9, 8, 1, 6, 9, 1, 6, 9, 1, 6] \text{ and } [8, 8, 1, 6, 9, 1, 6, 9, 1, 6]\)
		\end{enumerate}
		For each combination of $n$ and $m$, we run with 100 random starting points drawn from the constraints and then test the performance of BDCA vs DCA. \cref{fig: Set Scaling BDCA/DCA,fig: setscaling adaptive BDCA/DCA} show the results of the runs.
		
		As in \cref{Example: Scaling} both BDCA and adaptive BDCA are better than DCA in terms of iterations and run time. Similarly we observe the general trend of BDCA becoming increasingly faster compared to DCA as we increase the problem size and the evaluation of \(y_p\) in \cref{BDCA3} becomes more expensive. The set clustering problem is a more difficult problem than the basic clustering problem, and results in the BDCA being even more effective, particularly the adaptive BDCA. Notice that the average iteration and time ratios for both BDCA and adaptive BDCA in \cref{fig: Set Scaling BDCA/DCA,fig: setscaling adaptive BDCA/DCA} are nearly double of that in \cref{compare_BDCA,compare_BDCA_adap}. A table of the average runtimes and standard deviations can be found in \cref{table: DCA set,table: BDCA set,table: BDCA Adaptive set}. Overall, the results suggest that for set clustering problems at all scales, it is beneficial to use adaptive BDCA.
		
		\begin{figure}[htpb]
			\centering
			\begin{subfigure}[b]{0.49\textwidth}
				\centering
				\includegraphics[width=\textwidth]{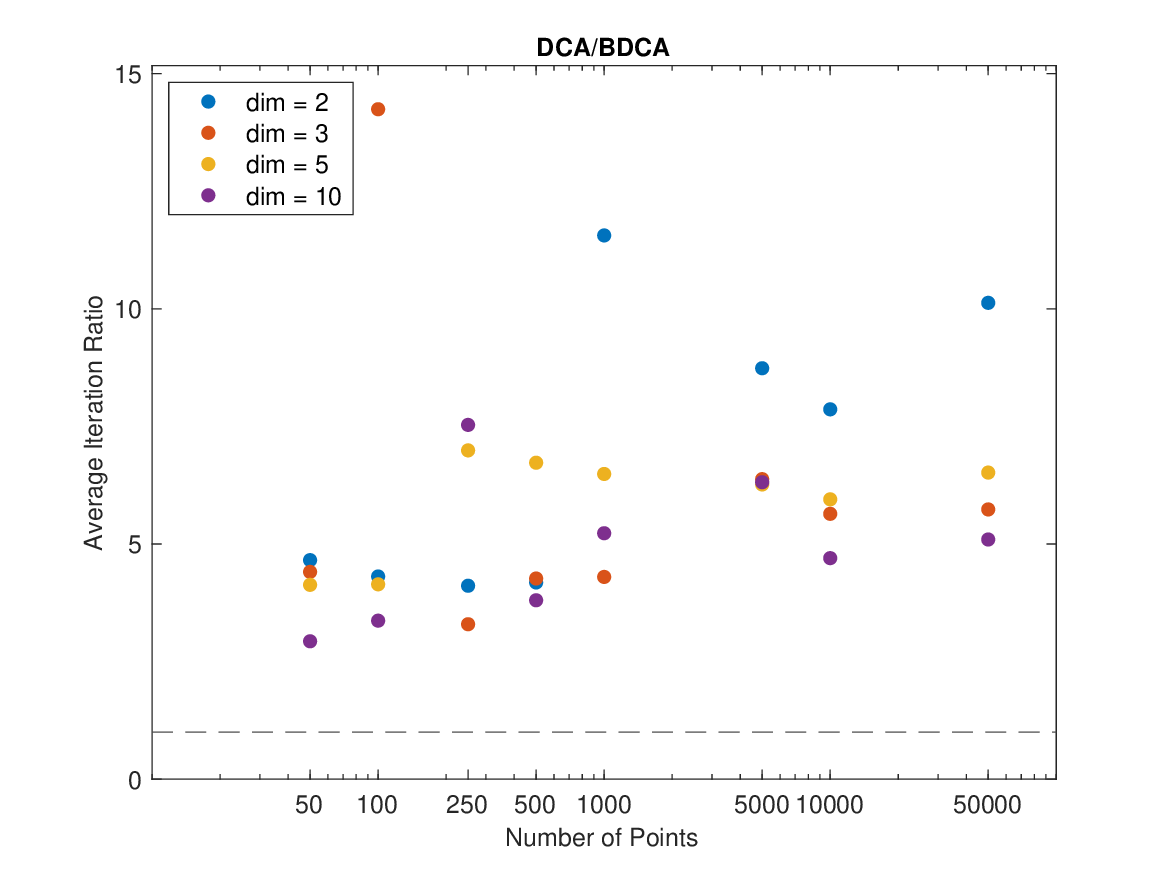}
			\end{subfigure}
			\hfill
			\begin{subfigure}[b]{0.49\textwidth}
				\centering
				\includegraphics[width=\textwidth]{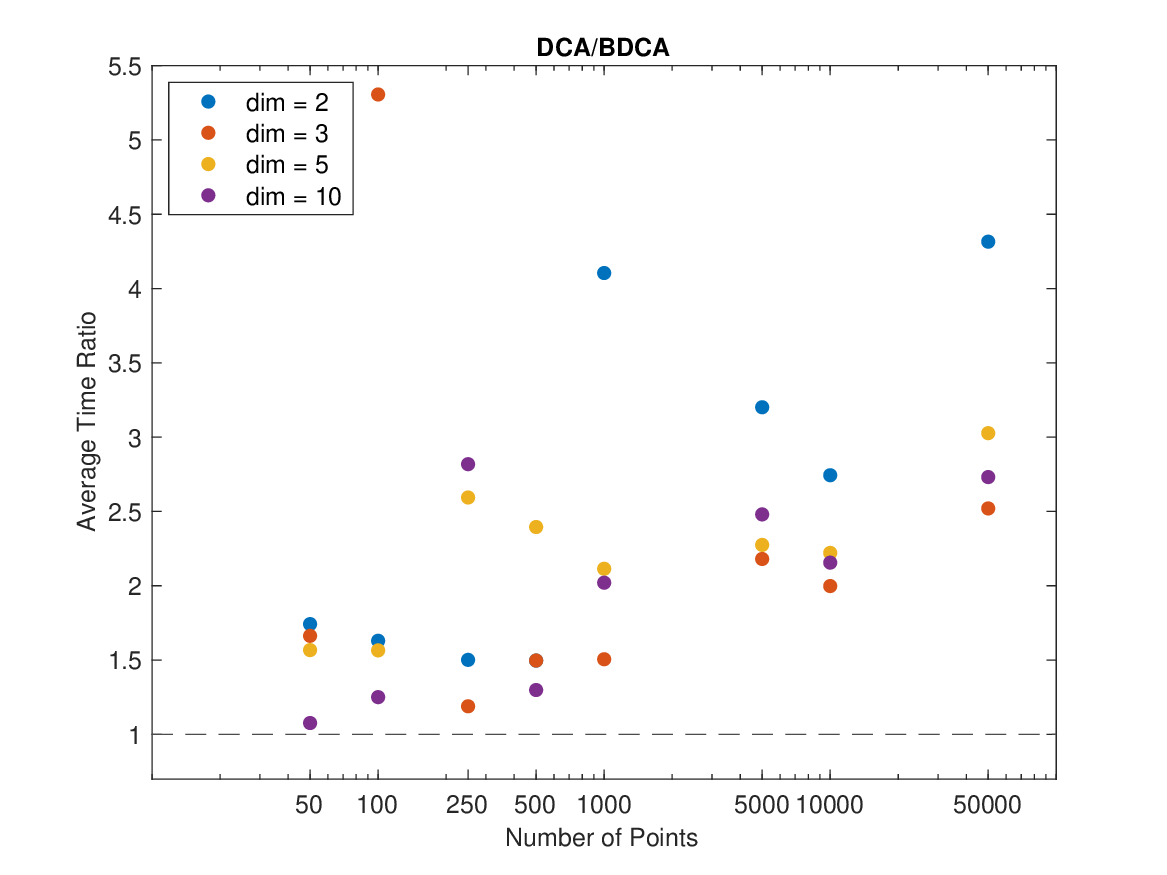}
			\end{subfigure}
			\caption{Iteration and Time Ratio Comparison between DCA and BDCA for \cref{Example: Set Clustering}.}
			\label{fig: Set Scaling BDCA/DCA}
		\end{figure}
		
		\begin{figure}[htpb]
			\centering
			\begin{subfigure}[b]{0.49\textwidth}
				\centering
				\includegraphics[width=\textwidth]{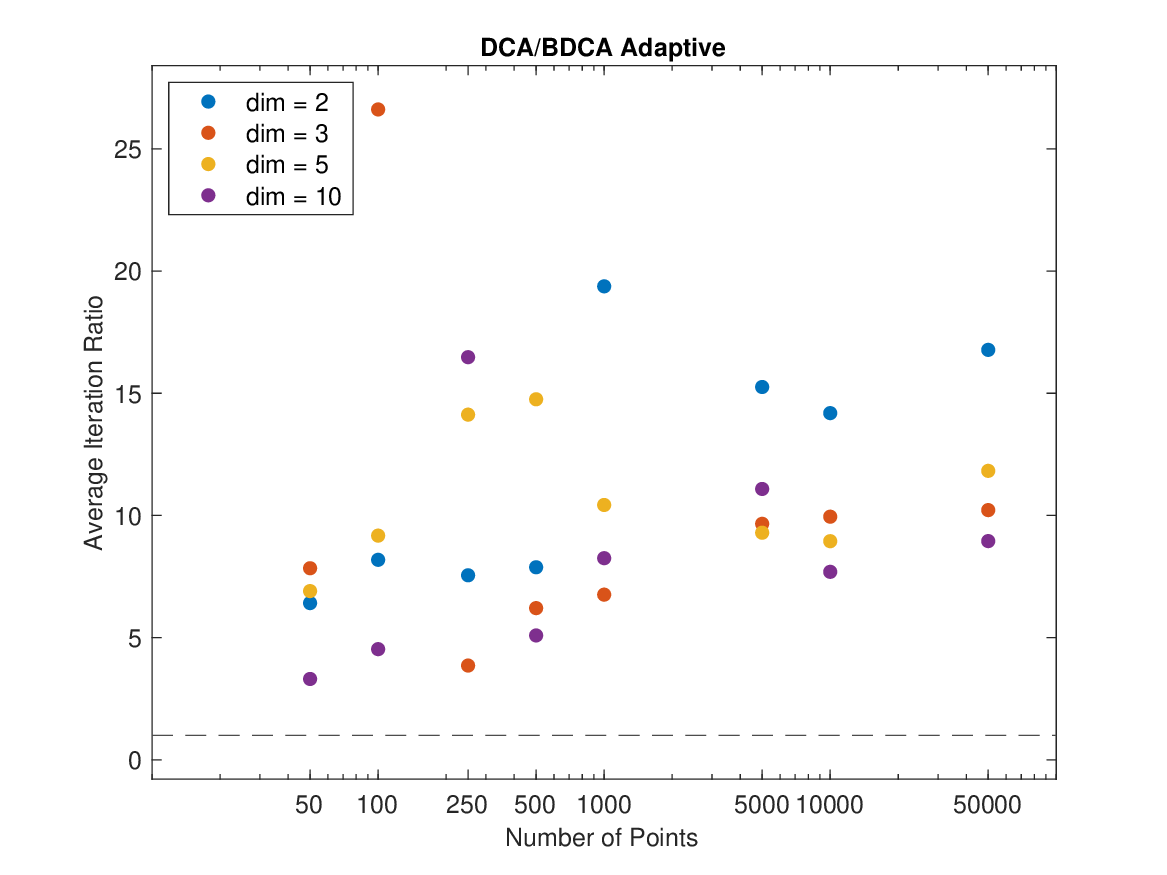}
			\end{subfigure}
			\hfill
			\begin{subfigure}[b]{0.49\textwidth}
				\centering
				\includegraphics[width=\textwidth]{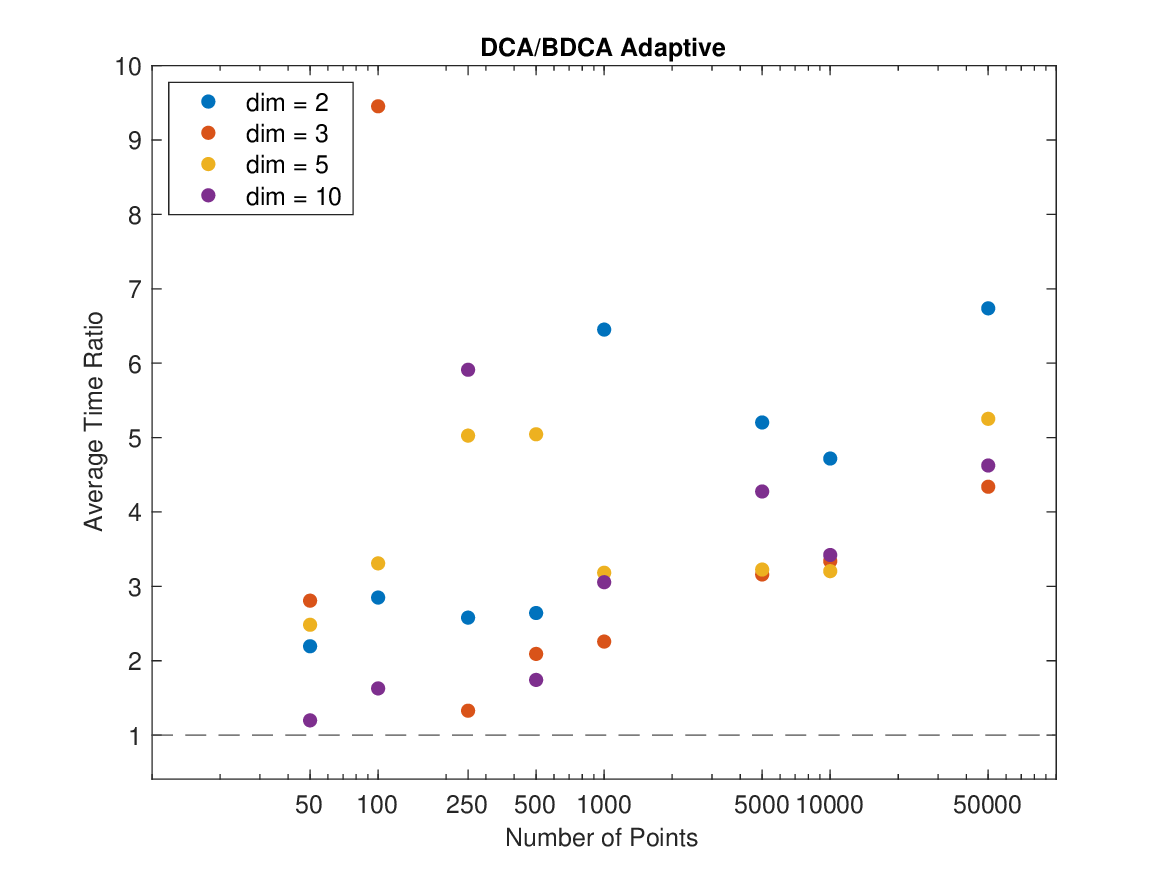}
			\end{subfigure}
			\caption{Iteration and Time Ratio Comparison between DCA and adaptive BDCA for \cref{Example: Set Clustering}.}
			\label{fig: setscaling adaptive BDCA/DCA}
		\end{figure}
		
		A visualization for \cref{Example: Set Clustering} when the number of points is $5000$ is in \cref{5000pts}.
		
		\begin{figure}[H]
			\centering
			\includegraphics[width=\textwidth]{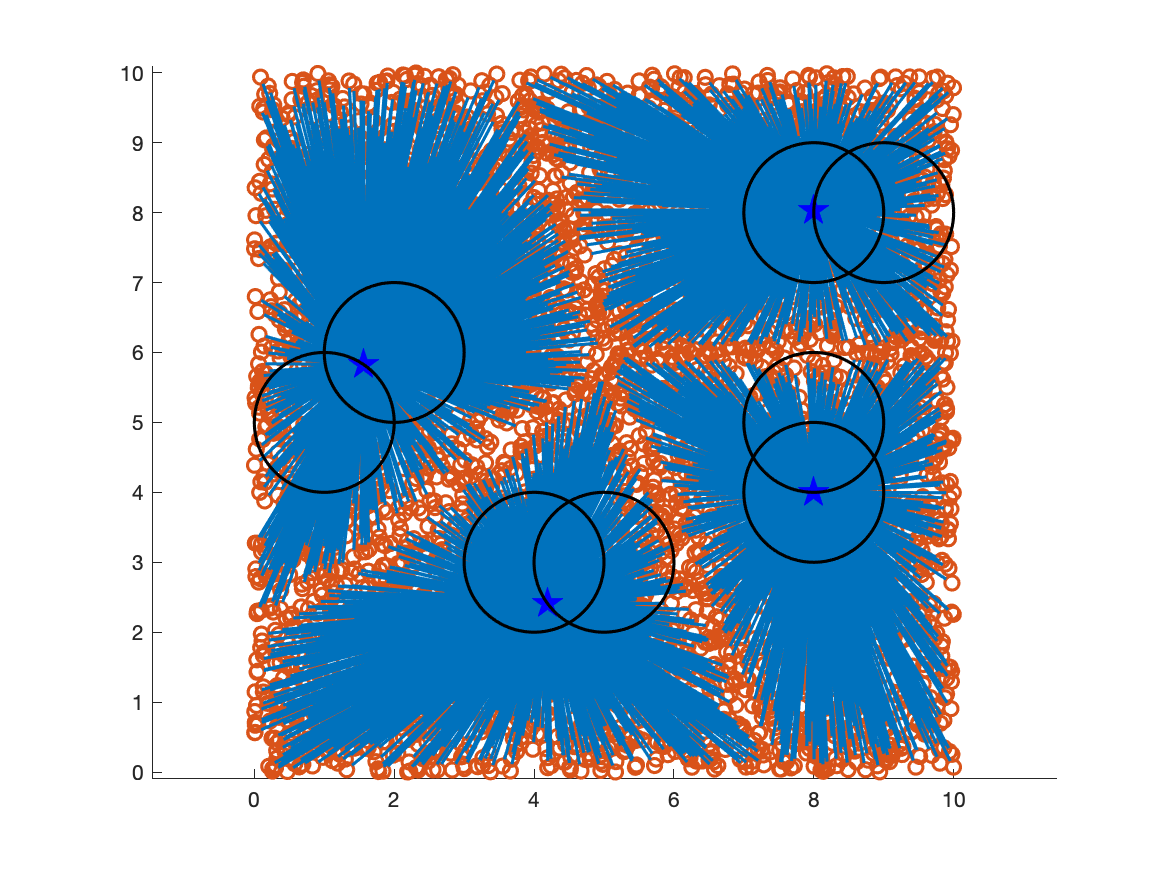}  
			\caption{A 4-center set clustering problems with 1500 points drawn from $[0,10]^2$. Each set is a ball with radius $0.1$, \cref{Example: Set Clustering}. Blue stars are centers.}
			\label{5000pts}
		\end{figure}
	\end{Example}	
	\section{Conclusions}\label{Chapter: Conclusions}
	The aim of this project was to investigate the application of the BDCA on clustering and set clustering problems with constraints. This is the first paper to test the application BDCA to any problem with nonlinear constraints. For each problem, we presented the DCA and penalty method used previously and suggested a BDCA-based method for solving it. We performed numerical experiments to test all of the methods described and presented our results in \cref{examples}, with the code and data from these examples available in the supplemental material. These experiments tested a variety of nonlinear constraints, and in all experiments, the BDCA method was able to achieve fewer iterations to convergence than the DCA method.  It also outperforms DCA in term of CPU running time.
	
	Overall, the work of this project has shown the potential effectiveness of BDCA-based methods for solving clustering and set clustering problems with constraints. The performance of these  algorithms is promising for application to practical clustering problems. Further experiments with higher dimensions and changing the number of  constraints could be another important direction for future work on this topic. Furthermore, investigations into accelerating the BDCA method, via problem dependent tuning of parameters is an area we will consider. 
	
  \bibliography{References/references.bib}

\begin{thebibliography}{10}
\expandafter\ifx\csname url\endcsname\relax
  \def\url#1{\texttt{#1}}\fi
\expandafter\ifx\csname urlprefix\endcsname\relax\def\urlprefix{URL }\fi
\expandafter\ifx\csname href\endcsname\relax
  \def\href#1#2{#2} \def\path#1{#1}\fi

\bibitem{An2005}
L.~T.~H. An, P.~D. Tao, \href{https://doi.org/10.1007/s10479-004-5022-1}{The {DC} ({Difference} of {Convex} {Functions}) {Programming} and {DCA} {Revisited} with {DC} {Models} of {Real} {World} {Nonconvex} {Optimization} {Problems}}, Annals of Operations Research 133~(1) (2005) 23--46.
\newblock \href {https://doi.org/10.1007/s10479-004-5022-1} {\path{doi:10.1007/s10479-004-5022-1}}.
\newline\urlprefix\url{https://doi.org/10.1007/s10479-004-5022-1}

\bibitem{TA1}
T.~P. Dinh, H.~A.~L. Thi, \href{https://www.semanticscholar.org/paper/Convex-analysis-approach-to-d.c.-programming%3A-and-Dinh-Thi/4e8861a2fc85ae0ad397ed385f8e00e18ac422e8}{Convex analysis approach to d.c. programming: {Theory}, {Algorithm} and {Applications}}, 1997.
\newline\urlprefix\url{https://www.semanticscholar.org/paper/Convex-analysis-approach-to-d.c.-programming%3A-and-Dinh-Thi/4e8861a2fc85ae0ad397ed385f8e00e18ac422e8}

\bibitem{LeThi2018}
H.~A. Le~Thi, T.~Pham~Dinh, {DC} programming and {DCA}: thirty years of developments, Mathematical Programming 169~(1) (2018) 5--68.

\bibitem{AragonArtacho2018}
F.~J. Aragón~Artacho, R.~M.~T. Fleming, P.~T. Vuong, \href{https://doi.org/10.1007/s10107-017-1180-1}{Accelerating the {DC} algorithm for smooth functions}, Mathematical Programming 169~(1) (2018) 95--118.
\newblock \href {https://doi.org/10.1007/s10107-017-1180-1} {\path{doi:10.1007/s10107-017-1180-1}}.
\newline\urlprefix\url{https://doi.org/10.1007/s10107-017-1180-1}

\bibitem{Artacho2020}
F.~J.~A. Artacho, P.~T. Vuong, The boosted difference of convex functions algorithm for nonsmooth functions, {SIAM} Journal on Optimization 30~(1) (2020) 980--1006.
\newblock \href {https://doi.org/10.1137/18m123339x} {\path{doi:10.1137/18m123339x}}.

\bibitem{An2007}
L.~T.~H. An, M.~T. Belghiti, P.~D. Tao, \href{https://doi.org/10.1007/s10898-006-9066-4}{A new efficient algorithm based on {DC} programming and {DCA} for clustering}, Journal of Global Optimization 37~(4) (2007) 593--608.
\newblock \href {https://doi.org/10.1007/s10898-006-9066-4} {\path{doi:10.1007/s10898-006-9066-4}}.
\newline\urlprefix\url{https://doi.org/10.1007/s10898-006-9066-4}

\bibitem{Nam2018}
N.~M. Nam, N.~T. An, S.~Reynolds, T.~Tran, \href{https://doi.org/10.1080/02331934.2018.1510498}{Clustering and multifacility location with constraints via distance function penalty methods and dc programming}, Optimization 67~(11) (2018) 1869--1894.
\newblock \href {https://doi.org/10.1080/02331934.2018.1510498} {\path{doi:10.1080/02331934.2018.1510498}}.
\newline\urlprefix\url{https://doi.org/10.1080/02331934.2018.1510498}

\bibitem{AragonArtacho2022}
F.~J. Aragón-Artacho, R.~Campoy, P.~T. Vuong, \href{https://doi.org/10.1007/s11228-022-00656-x}{The {Boosted} {DC} {Algorithm} for {Linearly} {Constrained} {DC} {Programming}}, Set-Valued and Variational Analysis 30~(4) (2022) 1265--1289.
\newblock \href {https://doi.org/10.1007/s11228-022-00656-x} {\path{doi:10.1007/s11228-022-00656-x}}.
\newline\urlprefix\url{https://doi.org/10.1007/s11228-022-00656-x}

\bibitem{Ferreira2021}
O.~P. Ferreira, E.~M. Santos, J.~C.~O. Souza, \href{http://arxiv.org/abs/2103.10757}{Boosted scaled subgradient method for {DC} programming}, Tech. rep., arXiv:2103.10757 [math] type: article (Mar. 2021).
\newblock \href {https://doi.org/10.48550/arXiv.2103.10757} {\path{doi:10.48550/arXiv.2103.10757}}.
\newline\urlprefix\url{http://arxiv.org/abs/2103.10757}

\bibitem{Ferreira2022}
O.~P. Ferreira, E.~M. Santos, J.~C.~O. Souza, \href{http://arxiv.org/abs/2111.01290}{A boosted {DC} algorithm for non-differentiable {DC} components with non-monotone line search}, Tech. rep., arXiv:2111.01290 [math] type: article (Jun. 2022).
\newblock \href {https://doi.org/10.48550/arXiv.2111.01290} {\path{doi:10.48550/arXiv.2111.01290}}.
\newline\urlprefix\url{http://arxiv.org/abs/2111.01290}

\bibitem{HUL}
J.-B. Hiriart-Urruty, C.~Lemar{\'{e}}chal, Fundamentals of Convex Analysis, Springer Berlin Heidelberg, 2001.
\newblock \href {https://doi.org/10.1007/978-3-642-56468-0} {\path{doi:10.1007/978-3-642-56468-0}}.

\bibitem{bn}
B.~S. Mordukhovich, N.~M. Nam, An Easy Path to Convex Analysis and Applications, Springer International Publishing, 2014.
\newblock \href {https://doi.org/10.1007/978-3-031-02406-1} {\path{doi:10.1007/978-3-031-02406-1}}.

\bibitem{r}
R.~T. Rockafellar, Convex {Analysis}, Princeton University Press, 1997, google-Books-ID: 1TiOka9bx3sC.

\bibitem{bmnt19}
A.~Bajaj, B.~S. Mordukhovich, N.~M. Nam, T.~Tran, \href{https://doi.org/10.1080/10556788.2020.1771335}{Solving a continuous multifacility location problem by {DC} algorithms}, Optimization Methods and Software 37~(1) (2022) 338--360.
\newblock \href {https://doi.org/10.1080/10556788.2020.1771335} {\path{doi:10.1080/10556788.2020.1771335}}.
\newline\urlprefix\url{https://doi.org/10.1080/10556788.2020.1771335}

\bibitem{TA2}
P.~D. Tao, L.~T.~H. An, \href{https://epubs.siam.org/doi/10.1137/S1052623494274313}{A {D}.{C}. {Optimization} {Algorithm} for {Solving} the {Trust}-{Region} {Subproblem}}, SIAM Journal on Optimization 8~(2) (1998) 476--505.
\newblock \href {https://doi.org/10.1137/S1052623494274313} {\path{doi:10.1137/S1052623494274313}}.
\newline\urlprefix\url{https://epubs.siam.org/doi/10.1137/S1052623494274313}

\bibitem{AnNam}
N.~T. An, N.~M. Nam, \href{https://doi.org/10.1080/02331934.2016.1253694}{Convergence analysis of a proximal point algorithm for minimizing differences of functions}, Optimization 66~(1) (2017) 129--147.
\newblock \href {https://doi.org/10.1080/02331934.2016.1253694} {\path{doi:10.1080/02331934.2016.1253694}}.
\newline\urlprefix\url{https://doi.org/10.1080/02331934.2016.1253694}

\bibitem{TAN}
H.~A. Le~Thi, V.~N. Huynh, T.~Pham~Dinh, \href{https://doi.org/10.1007/s10957-018-1345-y}{Convergence {Analysis} of {Difference}-of-{Convex} {Algorithm} with {Subanalytic} {Data}}, Journal of Optimization Theory and Applications 179~(1) (2018) 103--126.
\newblock \href {https://doi.org/10.1007/s10957-018-1345-y} {\path{doi:10.1007/s10957-018-1345-y}}.
\newline\urlprefix\url{https://doi.org/10.1007/s10957-018-1345-y}

\bibitem{Toland1979}
J.~F. Toland, \href{http://eudml.org/doc/94800}{On subdifferential calculus and duality in non-convex optimization}, Mémoires de la Société Mathématique de France 60 (1979) 177--183.
\newline\urlprefix\url{http://eudml.org/doc/94800}

\bibitem{CZL12}
E.~C. Chi, H.~Zhou, K.~Lange, \href{https://doi.org/10.1007/s10107-013-0697-1}{Distance majorization and its applications}, Mathematical Programming 146~(1) (2014) 409--436.
\newblock \href {https://doi.org/10.1007/s10107-013-0697-1} {\path{doi:10.1007/s10107-013-0697-1}}.
\newline\urlprefix\url{https://doi.org/10.1007/s10107-013-0697-1}

\bibitem{numopt}
S.~J.~W. Jorge~Nocedal, Numerical Optimization, Springer New York, 2006.
\newblock \href {https://doi.org/10.1007/978-0-387-40065-5} {\path{doi:10.1007/978-0-387-40065-5}}.

\bibitem{Nam2017}
N.~M. Nam, R.~B. Rector, D.~Giles, \href{https://doi.org/10.1007/s10957-017-1075-6}{Minimizing {Differences} of {Convex} {Functions} with {Applications} to {Facility} {Location} and {Clustering}}, Journal of Optimization Theory and Applications 173~(1) (2017) 255--278.
\newblock \href {https://doi.org/10.1007/s10957-017-1075-6} {\path{doi:10.1007/s10957-017-1075-6}}.
\newline\urlprefix\url{https://doi.org/10.1007/s10957-017-1075-6}

\bibitem{rei}
G.~Reinelt, {TSPLIB}{\textemdash}a traveling salesman problem library, {ORSA} Journal on Computing 3~(4) (1991) 376--384.
\newblock \href {https://doi.org/10.1287/ijoc.3.4.376} {\path{doi:10.1287/ijoc.3.4.376}}.

\end{thebibliography}
\end{document}